\documentclass[a4,semcolor,psfig,english,11pt,epsf,portrait]{article}
\usepackage{amsmath,amsfonts,latexsym,amscd,amssymb,theorem}

\usepackage[T1]{fontenc}
\usepackage{epsfig}
\usepackage{amsmath}
\def\R{\hbox{\bf R}}

\def\T{\mathcal{ T}}

\def\i{\infty}

\def\<{\langle}
\def\>{\rangle}
\newcommand{\ba}{\begin{eqnarray}}
\newcommand{\ea}{\end{eqnarray}}

\newtheorem{thm}{Theorem}[section]

\newtheorem{theorem}[thm]{Theorem}
\newtheorem{definition}[thm]{Definition}
\newtheorem{lemma}[thm]{Lemma}
\newtheorem{proposition}[thm]{Proposition}
\newtheorem{corollary}[thm]{Corollary}
\newtheorem{rem}[thm]{Remark}

\newcommand{\vr}{\varepsilon}
\newcommand{\eps}{\epsilon}
\numberwithin{equation}{section}

\renewcommand{\R}{{\mathbb R}}

\begin{document}

\title{\bf Existence and uniqueness for a nonlinear
  parabolic/Hamilton-Jacobi coupled system describing the dynamics of
  dislocation densities }
\author{
\normalsize\textsc{ Hassan
IBRAHIM\footnote{Cermics, Ecole des Ponts, ParisTech,
6 et 8 avenue Blaise Pascal, Cit\'e Descartes
 Champs-sur-Marne, 77455 Marne-la-Vall\'ee Cedex 2, France}
}}
\vspace{20pt}
\maketitle

%%%%%%%%%%%%%%%%%%%%%%%%%%%%%%%%%%%%%%%%%%%%%%%%%%%%%%%%%%%%%%%%%%%%%%%%%%
%%%%%%%%%%%%%%%%%%%%%%%%%%%%%%%%%%%%%%%%%%%%%%%%%%%%%%%%%%%%%%%%%%%%%%%%%%

 \centerline{\small{\bf{Abstract}}}
 \noindent{\small{We study a mathematical model describing the dynamics
of dislocation densities in crystals. This model is expressed as a
one-dimensional system of a parabolic equation and a first order Hamilton-Jacobi
equation that are coupled together. We show the existence and
uniqueness of a viscosity solution among those assuming a lower-bound on
their gradient for all time including the initial data. Moreover, we show the
existence of a viscosity solution when we have no such restriction on
the initial data. We also state a result of existence and uniqueness of an
entropy solution of the system obtained by spatial derivation. The 
uniqueness of this entropy solution holds in the class of ``bounded from
below'' solutions. In order to prove these results, we use a relation between scalar conservation laws and
Hamilton-Jacobi equations, mainly to get some gradient
estimates. This study will take place in $\R$, and on a bounded domain with
suitable boundary conditions.}}

\centerline{\small{\bf{Resumé}}}
 \noindent{\small{Nous étudions un modèle mathématique décrivant la dynamique de densités
de dislocations dans les cristaux. Ce modèle s'écrit comme un système
1D  couplant une équation parabolique et une équation de Hamilton-Jacobi
du premier ordre. On montre l'existence et l'unicité d'une solution de
viscosité dans la classe des fonctions ayant un gradient minoré pour
tout temps ainsi qu'au temps initial. De plus, on montre l'existence
d'une solution de viscosité sans cette condition sur la donnée initiale. On présente
également un résultat d'existence et d'unicité pour une solution
entropique d'un système obtenu par dérivation spatiale. L'unicité de
cette solution entropique a lieu dans la classe des solutions minorées.
Pour montrer ces résultats, on utilise une relation entre les lois de conservation scalaire et les
équations de Hamilton-Jacobi, principalement pour obtenir des contrôles
du gradient. Cette étude a lieu dans $\R$ et dans un domaine
borné avec des conditions aux bords appropriées.}}

\hfill\break
 \noindent{\small{\bf{AMS Classification: }}} {\small{70H20, 35L65,
49L25, 54C70, 74H20, 74H25.}}\hfill\break
  \noindent{\small{\bf{Key words: }}} {\small{Hamilton-Jacobi equations, scalar conservation laws, viscosity
solutions, entropy solutions, dynamics of dislocation densities.}}\hfill\break

%%%%%%%%%%%%%%%%%%%%%%%%%%%%%%%%%%%%%%%%%%%%%%%%%%%%%%%%%%%%%%%%%%%%%%%%%%
%%%%%%%%%%%%%%%%%%intoduction%%%%%%%%%%%%%%%%%%%%%%%%%%%%%%%%%%%%%%%%%%%
\begin{section}{Introduction}\label{sec1}
\subsection{Physical motivation}
A dislocation is a defect, or irregularity within a crystal
structure that can be observed by electron microscopy. The theory
was originally developed by Vito Volterra in 1905. Dislocations are
a non-stationary phenomena and their motion is the main explanation of the plastic deformation
in metallic crystals (see \cite{nab, HL} for a recent and mathematical presentation).

Geometrically, each dislocation is characterized by a physical
quantity called the Burgers vector, which is responsible for its
orientation and magnitude. Dislocations are classified as being
positive or negative due to the orientation of its Burgers vector,
and they can move in certain crystallographic directions.

Starting from the motion of individual dislocations, a continuum
description can be derived by adopting a formulation of dislocation
dynamics in terms of appropriately defined dislocation densities,
namely the density of positive and negative dislocations. In this
paper we are interested in the model described by Groma, Csikor and
Zaiser \cite{gcz}, that sheads light on the evolution of the
dynamics of the ``two type'' densities of a system of straight parallel
dislocations, taking into consideration the
influence of the short range dislocation-dislocation interactions.
The model was originally presented in $\mathbb{R}^{2}\times (0,T)$
as follows:
\begin{equation}\label{original}
\left\{
\begin{aligned}
\frac{\partial\theta^{+}}{\partial
  t}+\bold{b}\cdot\frac{\partial}{\partial\bold{r}}\left[\theta^{+}
  \left\{(\tau_{sc}+\tau_{eff})-AD\frac{\bold{b}}{(\theta^{+}+\theta^{-})}
  \cdot\frac{\partial}{\partial\bold{r}}\left(\theta^{+}-\theta^{-}\right)\right\}\right]=0,\\
\frac{\partial\theta^{-}}{\partial
  t}-\bold{b}\cdot\frac{\partial}{\partial\bold{r}}\left[\theta^{-}
  \left\{(\tau_{sc}+\tau_{eff})-AD\frac{\bold{b}}{(\theta^{+}+\theta^{-})}
  \cdot\frac{\partial}{\partial\bold{r}}\left(\theta^{+}-\theta^{-}\right)\right\}\right]=0.
\end{aligned}
\right.
\end{equation}
Where $T>0$,  $\bold{r}=(x,y)$ represents the spatial variable,
$\bold{b}$ is the burger's vector, $\theta^{+}(\bold{r},t)$ and
$\theta^{-}(\bold{r},t)$ denote the densities of the positive and
negative dislocations respectively. The quantity $A$ is defined by
the formula $A=\mu/[2\pi(1-\nu)]$, where $\mu$ is the shear modulus
and $\nu$ is the Poisson ratio.  $D$ is a non-dimensional constant.
Stress fields are represented through  the self-consistent stress
$\tau_{sc}(\bold{r},t)$, and  the effective stress
$\tau_{eff}(\bold{r},t)$. $\frac{\partial}{\partial\bold{r}}$
denotes the gradient with respect to the coordinate vector
$\bold{r}$. An earlier investigation of the continuum description
of the dynamics of dislocation densities has been done in
\cite{GB}. However, a major drawback of these investigations 
is that the short range dislocation-dislocation correlations have
been neglected and dislocation-dislocation interactions were
described only by the long-range term which is the self-consistent
stress field. Moreover, for the model  described in \cite{GB}, we
refer the reader to \cite{elhajj,elhfor} for a one-dimensional
mathematical and numerical study, and to \cite{CEMR} for a
two-dimensional existence result.

In our work, we are interested in a particular setting of
(\ref{original}) where we make the following
assumptions:
\begin{itemize}
\item[(a1)] the quantities in equations (\ref{original}) are independent of
$y$,
\item[(a2)] $\bold{b}=(1,0)$, and the constants $A$ and $D$ are set to be
$1$,
\item[(a3)] the effective stress is assumed to be zero.
\end{itemize}
\begin{rem} (a1) gives that the self-consistent stress $\tau_{sc}$
is null; this is a consequence of the definition of $\tau_{sc}$ (see
\cite{gcz}).
\end{rem}
Assumptions (a1)-(a2)-(a3) permit rewriting the original model as a
$\bold{1D}$ problem in $\mathbb{R}\times(0,T)$:
\begin{equation}\label{gcz}
\left\{
\begin{aligned}
&\theta^{+}_{t}(x,t)-\left(\theta^{+}(x,t)\left(\frac{\theta^{+}_{x}(x,t)-\theta^{-}_{x}(x,t)}{\theta^{+}(x,t)+\theta^{-}(x,t)}\right)\right)_{x}=0,\\
&\theta^{-}_{t}(x,t)+\left(\theta^{-}(x,t)\left(\frac{\theta^{+}_{x}(x,t)-\theta^{-}_{x}(x,t)}{\theta^{+}(x,t)+\theta^{-}(x,t)}\right)\right)_{x}=0.
\end{aligned}
\right.
\end{equation}
We consider an integrated form of (\ref{gcz}) and we let:
\begin{equation}\label{rhotheta}
\rho_{x}^{\pm}=\theta^{\pm},\quad\theta=\theta^{+}+\theta^{-},\quad\rho=\rho^{+}-\rho^{-}\;
\quad\mbox{and}\quad\kappa=\rho^{+}+\rho^{-},
\end{equation}
in order to obtain, for special values of the constants of integration, the
following system of PDEs in terms of $\rho$ and $\kappa$ :\\
\begin{equation}\label{mainhj}
\left\{
\begin{aligned}
&\kappa_{t}\kappa_{x}=\rho_{t}\rho_{x}\quad&\mbox{in}&\quad
Q_T=\mathbb{R}\times(0,T),\\
&\kappa(x,0)=\kappa^{0}(x)\quad&\mbox{in}&\quad \mathbb{R},
\end{aligned}
\right.
\end{equation}
and
\begin{equation}\label{mainheat}
\hspace{-2.2cm}\left\{
\begin{aligned}
&\rho_{t}=\rho_{xx} \quad&\mbox{in}&\quad Q_T,\\
&\rho(x,0)=\rho^{0}(x) \quad&\mbox{in}&\quad \mathbb{R},
\end{aligned}
\right.
\end{equation}\\
where $T>0$ is a  fixed constant. Enough regularity on the initial data
will be given in order to impose the physically relevant condition,\\
\begin{equation}\label{condition}
\kappa^{0}_{x}\geq |\rho^{0}_{x}|\;.
\end{equation}\
This condition is natural: it indicates nothing but the
positivity of the dislocation densities $\theta^{\pm}(x,0)$ at the
initial time (see (\ref{rhotheta})).

\subsection{Main results}
In this paper, we show the existence and uniqueness of a viscosity
solution $\kappa$ of (\ref{mainhj}) in the class of all Lipschitz continuousviscosity
solutions having special ``bounded from below'' spatial
gradients. However,
we show the existence of a Lipschitz continuousviscosity solution of
(\ref{mainhj}) when this restriction is relaxed. A relation
between scalar conservation laws and Hamilton-Jacobi equations will
be exploited to get almost all our gradient controls of $\kappa$.
This relation, that will be made precise later, will also lead to a
result of existence and uniqueness of a bounded entropy solution of
the following equation:
\begin{equation}\label{thetader} \left\{
\begin{aligned}
&\theta_{t}=\left(\frac{\rho_{x}\rho_{xx}}{\theta}\right)_{x}\quad&\mbox{in}&\quad
Q_T,\\
&\theta(x,0)=\theta^{0}(x)\quad&\mbox{in}&\quad\mathbb{R},
\end{aligned}
\right.
\end{equation}
which is deduced formally  by taking a spatial derivation of
(\ref{mainhj}). The uniqueness of this entropy solution is always restricted to the
class of bounded entropy solutions with a special lower-bound.

Let $Lip(\R)$ denotes:
$$Lip(\R)=\{f:\R\mapsto\R;\;\;f\;\mbox{is a  Lipschitz continuous
  function}\}.$$ 
We prove the following theorems:
\begin{theorem}{\bf (Existence and uniqueness of a viscosity
    solution)}\label{main1}\\
Let $T>0$. Take  $\kappa^{0}\in
Lip(\mathbb{R})$ and $\rho^{0}\in
C^{\infty}_{0}(\mathbb{R})$ as initial data that satisfy:
\begin{equation}\label{res_cond_ini_vis}
\kappa_{x}^{0}\geq \sqrt{(\rho_{x}^{0})^{2}+\epsilon^{2}}\quad \mbox{
a.e. in }\quad \R,
\end{equation}
for some constant $\eps>0$.
Then, given the solution $\rho$ of  (\ref{mainheat}),
there exists a viscosity solution $\kappa\in Lip(\bar{Q}_T)$ of
(\ref{mainhj}), unique among the viscosity solutions satisfying:
$$ \kappa_{x}\geq \sqrt{\rho_{x}^{2}+\epsilon^{2}}\quad\mbox{ a.e.  in }
\quad\bar{Q}_T.$$
\end{theorem}
\begin{theorem}{\bf(Existence and uniqueness of an entropy
    solution)}\label{main2}\\
Let $T>0$. Take $\theta^{0}\in L^{\infty}(\mathbb{R})$ and
$\rho^{0}\in C^{\infty}_{0}(\R)$ such that,
$$\theta^{0}\geq\sqrt{(\rho^{0}_{x})^{2}+\epsilon^{2}}\quad\mbox{
  a.e. in } \quad\R,$$
for some constant $\eps>0$. Then, there exists an entropy
solution $\theta\in L^{\infty}(\bar{Q}_T)$ of (\ref{thetader}), unique 
among the entropy solutions satisfying:
$$ \theta\geq \sqrt{\rho^{2}_{x}+\epsilon^{2}}\quad \mbox{
  a.e.}\;\mbox{ in }\quad \bar{Q}_T.$$
Moreover, we have $\theta=\kappa_x$, where $\kappa$ is the solution
given by Theorem \ref{main1}.
\end{theorem}
The notion of viscosity solutions and entropy solutions will be
recalled in Section \ref{sec2}. We now relate these results to our
one-dimensional problem (\ref{gcz}). Remarking that
$\rho_{x}=\theta^{+}-\theta^{-}$ and
$\kappa_{x}=\theta^{+}+\theta^{-}$, we have as a consequence:
\begin{corollary}{\bf(Existence and uniqueness for problem
(\ref{gcz}))}\label{reglast}\\
Let $T>0$. Let $\theta^{+}_{0}$ and $\theta^{-}_{0}$ be two given functions
representing the initial positive and negative dislocation densities
respectively. If the following conditions are satisfied:
\begin{itemize}
\item[(1)] $\theta^{+}_{0}-\theta^{-}_{0}\in C^{\infty}_{0}(\R),$
\item[(2)] $\theta^{+}_{0}$, $\theta^{-}_{0}\in L^{\infty}(\mathbb{R})$,
\end{itemize}
together with,
$$\theta^{+}_{0}+\theta^{-}_{0}\geq\sqrt{(\theta^{+}_{0}-\theta^{-}_{0})^{2}+\epsilon^{2}}\quad
\mbox{ a.e. in }\quad\R,$$ then there exists a solution
$(\theta^{+},\theta^{-})\in (L^{\infty}(Q_T))^{2}$ to the
system (\ref{gcz}), in the sense of Theorems \ref{main1} and
\ref{main2}, unique among those satisfying:
$$\theta^{+}+\theta^{-}\geq\sqrt{(\theta^{+}-\theta^{-})^{2}+\epsilon^{2}}\quad
\mbox{a.e.}\quad \mbox{in} \quad \bar{Q}_T.$$
\end{corollary}
\begin{rem}
Conditions (1) and (2) are sufficient requirements for the
compatibility with the regularity of $\rho^{0}$ and $\kappa^{0}$
previously stated.
\end{rem}
\begin{theorem}{\bf(Existence of a viscosity solution, case
$\eps=0$)}\label{maintheorem}\\
Let $T>0$, $\kappa^{0}\in Lip(\mathbb{R})$ and $\rho^{0}\in
C^{\infty}_{0}(\mathbb{R})$. If the condition
(\ref{condition}) is satisfied a.e. in $\R$, then there exists a
viscosity solution $\kappa\in Lip(\bar{Q}_T)$ of (\ref{mainhj})
satisfying:
\begin{equation}\label{for_all_time}
\kappa_{x}\geq|\rho_{x}|\quad\mbox{ a.e. in }\quad \bar{Q}_{T}.
\end{equation}
\end{theorem}
\begin{rem}
In the limit case where $\eps=0$, we remark that having
(\ref{for_all_time}) was intuitively expected due to the positivity
of the dislocation densities $\theta^{+}$ and $\theta^{-}$. This
reflects in some way the well-posedness of the model (\ref{gcz}) of
the dynamics of dislocation densities. We also remark that our
result of existence of a solution of (\ref{mainhj}) under
(\ref{for_all_time}) still holds if we start with
$\kappa^{0}_{x}=\rho^{0}_{x}=0$ on some interval of the real line.
In other words, we can imagine that we start with the probability of
the formation of no dislocation zones.
\end{rem}
\noindent \underline{\textbf{Problem with boundary conditions}}.\\
We consider once again problem (\ref{mainhj}), similar results to that announced above will be shown on a bounded
interval of the real line with Dirichlet boundary conditions (see
Section \ref{sec5}). This problem
corresponds physically to the study of the dynamics of dislocation
densities in a part of a material with the geometry of a slab (see \cite{gcz}).
\subsection{Organization of the paper} The paper is organized as
follows. In Section \ref{sec2}, we start by stating the definition
of viscosity and entropy solutions with some of their properties. In
Section \ref{sec3}, we prove the existence and uniqueness of a
viscosity solution to an approximated problem of (\ref{mainhj}),
namely Proposition \ref{approblem}, and we move
on, giving additional properties of our approximated solution
(Proposition \ref{addpro}) and consequently proving Theorems
\ref{main1} and \ref{main2}. In Section \ref{sec4}, we present the
proof of Theorem \ref{maintheorem}. Section \ref{sec5} is devoted to
the study of problem (\ref{mainhj}) on a bounded domain with
suitable boundary conditions. Finally, Section \ref{sec6} is an
appendix containing a sketch of the proof to the classical
comparison principle of scalar conservation laws adapted to our
equation with low regularity.
\end{section}

%%%%%%%%%%%%%%%%%%%%%%%%%Notations and preliminaries%%%%%%%%%%%%%%%%%%%%%%%%%%%%%
%%%%%%%%%%%%%%%%%%%%%%%%%%%%%%%%%%%%%%%%%%%%%%%%%%%%%%%%%%%%%%%%%
\begin{section}{Notations and Preliminaries}\label{sec2} We first fix some
  notations. If $\Omega$ is an open subset of $\mathbb{R}^{n}$, $k$ is a
  positive integer, we denote by $C^{k}(\Omega)$ the
  space of all real valued $k$ times continuously differentiable
  functions. $C^{k}_{0}(\Omega)$ is the subspace of $C^{k}(\Omega)$
  consisting of function of compact support in $\Omega$, and
  $C^{k}_{b}(\Omega)=C^{k}(\Omega)\cap W^{k,\infty}(\Omega)$ where 
$W^{k,\infty}(\Omega)$ is
  defined below. Furthermore, let $UC(\Omega)$ and $Lip(\Omega)$ denote
  the spaces  of uniformly  continuous functions and  Lipschitz continuousfunctions on 
$\Omega$
  respectively. The sobolev space $W^{m,p}(\Omega)$ with $m\geq 1$ an
  integer and $p: 1\leq p \leq \infty$ a real,
  is defined by
\begin{eqnarray*}
W^{n,p}(\Omega)=\left\{u\in L^{p}(\Omega)\left|\begin{aligned}
\forall \alpha \mbox{ with } |\alpha|\leq n\;\; \exists f_{\alpha}\in
L^{p}(\Omega)\;\;\mbox{ such that }\\
\int_{\Omega}uD^{\alpha}\phi=(-1)^{|\alpha|}\int_{\Omega}f_{\alpha}\phi\;\;
\forall\phi\in C^{\infty}_{0}(\Omega)\end{aligned}\right.\right\},
\end{eqnarray*}
where we denote $D^{\alpha}u=f_{\alpha}$. This space
equipped with the norm
$$ ||u||_{W^{n,p}}=\sum_{0\leq |\alpha|\leq n}||D^{\alpha}u||_{L^{p}}$$
is a Banach space. In what follows, $T>0$. A map
$m:[0,\infty)\mapsto[0,\infty)$ that satisfy
\begin{itemize}
\item[$\bullet$] $m$ is continuous and non-decreasing;
\item[$\bullet$] $\displaystyle\lim_{x\rightarrow 0^{+}}m(x)=0$;
\item[$\bullet$] $m(a+b)\leq m(a)+m(b)$ for $a, b\geq 0$;
\end{itemize}
is said to be ``a modulus'', and $UC_{x}(\Omega\times[0,T])$ denotes
the space of those $u\in C(\Omega\times[0,T])$ for which there is a
modulus $m$ and $r>0$ such that $$ |u(x,t)-u(y,t)|\leq
m(|x-y|)\mbox{ for } x, y\in \Omega,\,|x-y|\leq r\mbox{ and }
t\in[0,T].$$ 

%%%%%%%%%%%%%%%%%%%%%%%%%%%%%%%%%%%%%%%%%%%%%%%%%%%%%%%%%%%%%%%%%%%%%%%%%%%%%%%%%%%%%%%%%%%%%%%%%

\noindent We will deal with two types of equations:\\

1. Hamilton-Jacobi equation:
\begin{equation}\label{HJ1}
\left\{
\begin{aligned}
&u_t+F(x,t,u_x)=0\quad&\mbox{in}&\quad Q_T,\\
&u(x,0)=u^{0}(x)\quad&\mbox{in}&\quad \mathbb{R},
\end{aligned}
\right.
\end{equation}

2. Scalar conservation laws:
\begin{equation}\label{CL}
\left\{
\begin{aligned}
&v_t+(F(x,t,v))_x=0\quad&\mbox{in}&\quad Q_T,\\
&v(x,0)=v^{0}(x)\quad&\mbox{in}&\quad \mathbb{R},
\end{aligned}
\right.
\end{equation}
where
$$
\begin{aligned}
F: \R\times[0&,T]\times\R\quad\rightarrow\quad\hspace{0.5cm} \R\\
(x,\,&t,\,u)\quad\quad\hspace{0.2cm}\mapsto\quad F(x,t,u)
\end{aligned}
$$
is called the Hamiltonian in the Hamilton-Jacobi equations and the
flux function in the scalar conservation laws. We will agree on the
continuity of this function, while additional and specific
regularity will be given when it is needed.
\begin{rem} We will use
the function $F$ as a notation for the Hamiltonian/flux function.
Although $F$ might differ from one equation to another, it will be
clarified in all what follows.
\end{rem}
\begin{rem}\label{7jinan}
The major part of this work concerns a Hamiltonian/flux function of a
special
form, namely:
\begin{equation}\label{shape_hamil}
F(x,t,u)=g(x,t)f(u),
\end{equation}
where such forms often arise in problems of physical interest
including traffic flow \cite{dantra} and two-phase flow in porous
media \cite{danpor}.
\end{rem}

%%%%%%%%%%%%%%%%%%%%%%%%%%%%%%%%%%%%%%%%%%%%%%%%%%%%%%%%%%%%%%%%%%%%%%%%%%%%%%%%

We start by defining the notion of viscosity solution to
Hamilton-Jacobi equations (\ref{HJ1}), and entropy solution to 
scalar
conservation laws (\ref{CL}) with a flux function given by Remark 
\ref{7jinan}, as well as some
results about existence, uniqueness, and regularity properties of
these solutions. We will end by a classical relation between these two
problems. These results will be needed throughout this
paper, precise references for the proofs will be mentioned later on.
\subsection{Viscosity solution: definition and properties}
\begin{definition}{\bf(\cite{CL2}, Viscosity solution: non-stationary 
case)}\\
1) A function $u\in C(Q_T;\mathbb{R})$ is a viscosity sub-solution of 
\begin{equation}\label{fir_of_HJ1}
u_t+F(x,t,u_x)=0\quad\mbox{in}\quad Q_T,
\end{equation}
if for every $\phi\in C^{1}(Q_T)$, whenever
$u-\phi$ attains a local maximum at $(x_0,t_0)\in Q_T$, then
$$\phi_t(x_0,t_0)+F(x_0,t_0,\phi_x(x_0,t_0))\leq 0.$$
2) A function $u\in C(Q_T;\mathbb{R})$ is a viscosity super-solution of 
(\ref{fir_of_HJ1}) if for every $\phi\in C^{1}(Q_T)$, whenever
$u-\phi$ attains a local minimum at $(x_0,t_0)\in Q_T$, then
$$\phi_t(x_0,t_0)+F(x_0,t_0,\phi_x(x_0,t_0))\geq 0.$$
3) A function $u\in C(Q_T;\mathbb{R})$ is a viscosity solution of
(\ref{fir_of_HJ1}) if it is both a viscosity sub- and super-solution of
(\ref{fir_of_HJ1}).\\
4) A function  $u\in C(\bar{Q}_T;\mathbb{R})$ is a viscosity solution of
the initial value problem (\ref{HJ1}) if $u$ is a viscosity
solution of (\ref{fir_of_HJ1}) and $u(x,0)=u^0(x)$ in $\mathbb{R}$.
\end{definition}
It is worth mentioning here that if a viscosity solution of a
Hamilton-Jacobi equation is differentiable at a certain point, then it
solves the equation there (see \cite[Corollary I.6]{CL2}). An equivalent
definition depending on the sub- and
super-differential of a continuous function is now presented. This
definition will be used for the demonstration of Proposition \ref{lipreghj}. 
Let us recall that the sub- and
the super-differential of a continuous function $u\in
C(\mathbb{R}^{n}\times (0,T))$, at a point $(x,t) \in
\mathbb{R}^{n}\times (0,T)$, are defined as the closed convex sets:\\
\begin{equation}\hspace{-6cm}D^{1,-}u(x,t)=\Big\{(p,\alpha)\in\mathbb{R}^{n}\times\mathbb{R}:
\nonumber
\end{equation}
\begin{equation}\hspace{2.6cm}\liminf_{(y,s)\rightarrow(x,t)}\frac{u(y,s)-u(x,t)-(p\cdot(y-x)+\alpha\cdot(s-t))}
{|y-x|+|s-t|}\geq 0\Big\}, \nonumber
\end{equation}
and\\
\begin{equation}\hspace{-6cm}D^{1,+}u(x,t)=\Big\{(p,\alpha)\in\mathbb{R}^{n}\times\mathbb{R}:
\nonumber
\end{equation}
\begin{equation}\hspace{2.6cm}\limsup_{(y,s)\rightarrow(x,t)}\frac{u(y,s)-u(x,t)-(p\cdot(y-x)+\alpha\cdot(s-t))}
{|y-x|+|s-t|}\leq 0\Big\}, \nonumber
\end{equation}
respectively.
\begin{definition}{\bf(Equivalent definition of viscosity solution)}\\
1) A function $u\in C(\mathbb{R}^{n}\times (0,T))$ is a viscosity
super-solution of (\ref{HJ1}) if and only if, for every $(x,t)\in
\mathbb{R}^{n}\times (0,T)$:
\begin{equation}\label{sousd}
\forall (p,\alpha)\in D^{1,-}u(x,t),\qquad \alpha+F(x,t,p)\geq 0.
\end{equation}
2) A function $u\in C(\mathbb{R}^{n}\times (0,T))$ is a viscosity
sub-solution of (\ref{HJ1}) if and only if, for every $(x,t)\in
\mathbb{R}^{n}\times (0,T)$:
\begin{equation}\label{surd}
\forall (p,\alpha)\in D^{1,+}u(x,t),\qquad \alpha+F(x,t,p)\leq 0.
\end{equation}
\end{definition}
\noindent This definition is more local, for it permits
verification that a given explicit function is a viscosity solution
in a more classical way, i.e. using the derivative calculus. A
similar definition, that will be used later, could be given in the
stationary case. Let $\Omega\subset\mathbb{R}^{n}$ be an open
domain, and consider the PDE
\begin{equation}\label{one}
F(x,u(x),\nabla u(x))=0, \quad \forall x\in \Omega,
\end{equation}
where
$F:\Omega\times\mathbb{R}\times\mathbb{R}^{n}\mapsto\mathbb{R}$ is a
continuous mapping.
\begin{definition}{\bf(Viscosity solution: stationary case)}\\
A continuous function $u\,:\,\Omega\mapsto\mathbb{R}$ is a
viscosity sub-solution of the PDE (\ref{one}) if for any continuously
differentiable function $\phi\,:\,\Omega\mapsto\mathbb{R}$ and any local
maximum $x_0\in\Omega$ of $u-\phi$, one has
$$F(x_{0},u(x_{0}),\nabla\phi(x_{0}))\leq 0.$$
Similarly, if at any local minimum point $x_0\in\Omega$ of $u-\phi$, one
has
$$ F(x_{0},u(x_{0}),\nabla\phi(x_{0}))\geq 0,$$
then $u$ is a viscosity super-solution.
Finally, if $u$ is both a viscosity sub-solution and a viscosity
super-solution, then $u$ is called a viscosity solution.
\end{definition}
%%%%%%%%%%%%%%%%%%%%%%%%%%%%%%%%%%%%%%%%%%%%%%%%%%%%%%%%%%%%%%%%%%%%%
In fact, this definition is used for interpreting solutions of
(\ref{mainhj}) in the viscosity sense. Furthermore, we say that $u$ is a
viscosity solution of the Dirichlet problem (\ref{one}) with
$u=\zeta\in C(\partial\Omega)$ if:\\
 (1) $u\in C(\bar{\Omega})$,\\
 (2) $u$ is a viscosity solution of (\ref{one}) in $\Omega$,\\
 (3) $u=\zeta$ on $\partial\Omega$.\\
For a better understanding of the
viscosity interpretation of boundary conditions of Hamilton-Jacobi
equations, we refer the reader
to \cite[Section 4.2]{barbook}.

Now, we will proceed by giving the main results concerning viscosity
solutions of (\ref{HJ1}). In order to have existence and
uniqueness, the Hamiltonian $F$ will be restricted by the following
conditions :\\\\
$\bold{(F0)}$ $F\in C(\mathbb{R}\times[0,T]\times\mathbb{R})$;\\\\
$\bold{(F1)}$ for each $R>0$ there is a constant $C_R$ such that for
all $(x,t,p)$, $(y,t,q)\in\mathbb{R}\times[0,T]\times[-R,R],$
$$ |\,F(x,t,p)-F(y,t,q)\,| \leq C_R(\,|p-q|+|x-y|);$$
$\bold{(F2)}$ there is a constant $C_F$ such that for all
$(t,p)\in[0,T]\times\mathbb{R}$ and all $x,\,y\in\mathbb{R},$
$$ |\,F(x,t,p)-F(y,t,p)\,|\leq C_F|x-y|(1+|p|).$$
We use these conditions to write down some results on viscosity
solutions.
\begin{theorem}\label{comphj}{\bf(Comparison, \cite[Theorem
1]{CL1})}\\
Let $F$ satisfy
$\bold{(F0)}$-$\bold{(F1)}$-$\bold{(F2)}$. If $u$, $\bar u\in
  UC_{x}(\bar{Q}_T)$ are two viscosity sub- and super-solution of
  the Hamilton-Jacobi equation (\ref{HJ1})
  respectively, with
$$ u(x,0)\leq\bar{u}(x,0)\quad \mbox{ in }\quad \mathbb{R},$$
then $u\leq\bar{u}$ in $\bar{Q}_T$.
\end{theorem}
\begin{theorem}\label{exishj}{\bf(Existence, \cite[Theorem 1]{CL1})}\\Let 
$F$ satisfy
$\bold{(F0)}$-$\bold{(F1)}$-$\bold{(F2)}$. If $u^0\in
  UC(\mathbb{R})$, then (\ref{HJ1}) has a viscosity solution 
$u\in
  UC_x(\bar{Q}_T)$.
\end{theorem}
\begin{rem}
The ``comparison'' theorem stated above gives the  uniqueness of the
viscosity solution.
\end{rem}
\begin{rem}\label{cond_V}
In the case where the Hamiltonian has the form
$$F(x,t,u)=g(x,t)f(u),$$
the following conditions:\\\\
$\bold{(V0)}$ $f\in C_{b}^{1}(\R;\R),$\\\\
$\bold{(V1)}$ $g \in C_{b}(\bar{Q}_T;\R),$\\\\
$\bold{(V2)}$ $g_{x}\in L^{\infty}(\bar{Q}_{T})$,\\\\
imply $\bold{(F0)}$-$\bold{(F1)}$-$\bold{(F2)}$ together with the boundedness
of the Hamiltonian.
\end{rem}
The next proposition reflects the behavior of viscosity
solutions under additional regularity assumptions on $u^0$ and $F$.
\begin{proposition}\label{lipreghj}{\bf(Additional regularity of the
  viscosity solution)}\\
Let $F=gf$ satisfy $\bold{(V0)}$-$\bold{(V1)}$-$\bold{(V2)}$. If
$u^0\in Lip(\mathbb{R})$ and $u\in UC_x(\bar{Q}_T)$
is the unique viscosity solution of (\ref{HJ1}), then
$u\in Lip(\bar{Q}_T)$.
\end{proposition}
\noindent {\bf Proof.} Consider the function $u^{\vr}$
defined on $\mathbb{R}\times [0,T]$ by:
$$u^{\vr}(x,t)=\sup_{y\in 
\R}\left\{u(y,t)-e^{kt}\frac{|x-y|^{2}}{2\vr}\right\}.$$
By \cite[Theorem 3]{ish}, the function $u$ satisfies,
$$|u(x,t)|\leq c^{*}(|x|+1)\quad \mbox{for}\;(x,t)\in \R\times[0,T],$$
where $c$ and $c^{*}$ are two positive constants. Therefore, $u$ is
a sublinear function for every time $t\in[0,T]$. The function
$u^{\vr}$ is defined via a supremum which is attained because of the
sublinearity of the function $u$ (a quadratic function always
control a linear one); the supremum can be achieved at several
points; let $x_{\vr}$ be one of them, so we can write
$$u^{\vr}(x,t)=u(x_{\vr},t)-e^{kt}\frac{|x-x_{\vr}|^{2}}{2\vr}.$$
We are going to prove that for $(p,\alpha)\in \R\times\R$, we have:
\begin{equation}\label{statement}
(p,\alpha)\in
D^{1,+}u^{\vr}(x,t)\,\Rightarrow\,\left(p,\alpha+ke^{kt}\frac{|x-x_{\vr}|^{2}}{2\vr}
\right)\in D^{1,+}u(x_{\vr},t).
\end{equation}
Since $(p,\alpha)\in D^{1,+}u^{\vr}(x,t)$, then we can write for
$(y,s)\sim(x,t)$ that,
\begin{equation}\label{ineq}
L=u^{\vr}(y,s)\leq u^{\vr}(x,t)+\alpha(s-t)+p(y-x)+o(|s-t|+|y-x|)=R,
\end{equation}
where the left side $L$ of (\ref{ineq}) satisfies,
\begin{equation}\label{ineq1}
L \geq u(z,s)-e^{ks}\frac{|z-y|^{2}}{2\vr},\quad z\in \R,
\end{equation}
and the right side $R$ of (\ref{ineq}) satisfies,
\begin{equation}\label{ineq2}
R\leq
u(x_{\vr},t)-e^{kt}\frac{|x-x_{\vr}|^{2}}{2\vr}+\alpha(s-t)+p(y-x)+o(|s-t|+|y-x|).
\end{equation}
Choose $z$ such that $z-y=x_{\vr}-x$, then
\begin{equation}\label{ineq3} z=x_{\vr}+(y-x)\sim
x_{\vr},\,\,\mbox{since}\,\,\, y\sim x.
\end{equation}
Combining (\ref{ineq}), (\ref{ineq1}), (\ref{ineq2}) and
(\ref{ineq3}) together, we get
\begin{equation}
\hspace{-6.8cm}u(x_{\vr}+(y-x),s)-e^{ks}\frac{|x-x_{\vr}|^{2}}{2\vr}\leq
\nonumber
\end{equation}
\begin{equation}
\hspace{1.4cm}u(x_{\vr},t)-e^{kt}\frac{|x-x_{\vr}|^{2}}{2\vr}+\alpha(s-t)+p(z-x_{\vr})+o(|s-t|+|z-x_{\vr}|),\nonumber
\end{equation}
and hence,
\begin{equation}
\hspace{-6cm}u(z,s)\leq
u(x_{\vr},t)+(e^{ks}-e^{kt})\frac{|x-x_{\vr}|^{2}}{2\vr}\nonumber
\end{equation}
\begin{equation}\label{ineq4}
\hspace{1.3cm}+\alpha(s-t)+p(z-x_{\vr})+o(|s-t|+|z-x_{\vr}|).
\end{equation}
We have
$$(e^{ks}-e^{kt})\frac{|x-x_{\vr}|^{2}}{2\vr}=ke^{kt}\frac{|x-x_{\vr}|^{2}}{2\vr}(s-t)+o(|s-t|),$$
then using inequality (\ref{ineq4}), we get
\begin{equation}
\hspace{-4.7cm}u(z,s)\leq
u(x_{\vr},t)+\left(\alpha+ke^{kt}\frac{|x-x_{\vr}|^{2}}{2\vr}\right)(s-t)\nonumber
\end{equation}
\begin{equation}
\hspace{3.301cm}+p(z-x_{\vr})+o(|s-t|+|z-x_{\vr}|),\nonumber
\end{equation}
which proves that $$\left(\alpha+ke^{kt}\frac{|x-x_{\vr}|^{2}}{2\vr}
,p\right)\in D^{1,+}u(x_{\vr},t),$$
and hence statement (\ref{statement}) is true.
Since $u$ is a viscosity sub-solution of (\ref{HJ1}), we have
$$\alpha+ke^{kt}\frac{|x-x_{\vr}|^{2}}{2\vr}+F(x_{\vr},t,p)\leq 0.$$
We use condition $(\bold{F1})$ with $p=q$, to get
\begin{eqnarray*}
\alpha+ke^{kt}\frac{|x-x_{\vr}|^{2}}{2\vr}+F(x,t,p)&\leq&F(x,t,p)-F(x_{\vr},t,p),\\
&\leq& C|x-x_{\vr}|,
\end{eqnarray*}
therefore,
\begin{eqnarray*}
\alpha+F(x,t,p)&\leq&
C|x-x_{\vr}|-ke^{kt}\frac{|x-x_{\vr}|^{2}}{2\vr},\\
&\leq& Cr_{\vr}-k\frac{r_{\vr}^{2}}{2\vr},\\
&\leq& \sup_{r>0}\left(Cr-\frac{kr^{2}}{2\vr}\right),
\end{eqnarray*}
where $r_{\vr}=|x-x_{\vr}|.$ At the maximum $\bar{r}$, we have
$C=\frac{k\bar{r}}{\vr}$. By choosing $k=\frac{C^2}{2}$, we get
$$\alpha+F(x,t,p)\leq \vr.$$
This inequality shows that $v^{\vr}=u^{\vr}-\vr t$ is a viscosity
sub-solution of (\ref{HJ1}) with
$v^{\vr}(x,0)=u^{\vr}(x,0)$. By the comparison principle, we have
\begin{eqnarray*}
v^{\vr}(x,t)-u(x,t)&\leq& \sup_{x\in
\R}(v^{\vr}(x,0)-u^{0}(x)),\\
&\leq& \sup_{x\in \R}(u^{\vr}(x,0)-u^{0}(x)),\\
&\leq& \sup_{x\in\R}\left(\sup_{y\in
\R}\left\{u^{0}(y)-\frac{|x-y|^{2}}{2\vr}\right\}-u^{0}(x)\right),\\
&\leq& \sup_{x,y\in \R}\left(\gamma|x-y|-\frac{|x-y|^{2}}{2\vr}\right),\\
&\leq& \sup_{r\geq0}\left(\gamma r-\frac{r^{2}}{2\vr}
\right)=\frac{\gamma^{2}\vr}{2},
\end{eqnarray*}
where $\gamma$ is the Lipschitz constant of the function $u^{0}$,
and $r=|x-y|$. This altogether shows the following inequality for
$x,y\in \R$:
\begin{equation}\label{jano}
u(y,t)-e^{kt}\frac{|x-y|^{2}}{2\vr}\leq u^{\vr}(x,t)\leq u(x,t)+\vr
t + \frac{\gamma^{2}\vr}{2}.
\end{equation}
Remark here that $k$ is a fixed; previously chosen constant.
Inequality (\ref{jano}) yields:
\begin{equation}\label{jojo}
u(y,t)-u(x,t)\leq e^{kt}\frac{|x-y|^{2}}{2\vr}+
\left(t+\frac{\gamma^{2}}{2} \right)\vr=\zeta/\vr + \beta\vr,
\end{equation}
where $\zeta=e^{kt}\frac{|x-y|^{2}}{2}$ and
$\beta=\left(t+\frac{\gamma^{2}}{2} \right)$. We minimize inequality
(\ref{jojo}) over $\vr$ to obtain,
\begin{eqnarray*}
u(y,t)-u(x,t)&\leq&
2\sqrt{\zeta\beta},\\
&\leq& e^{\frac{kt}{2}}\sqrt{2}\sqrt{t+\frac{\gamma^{2}}{2}}|x-y|.
\end{eqnarray*}
Since this inequality holds $\forall x, y\in \R$, exchanging $x$ with
$y$ yields,
$$|u(x,t)-u(y,t)|\leq C(F,u_0)|x-y| \quad \forall x, y \in
\R\;\;\mbox{and} \;\;t\in [0,T].$$
This shows that the function $u$ is Lipschitz continuous in $x$, uniformly in time
$t$. To prove the
Lipschitz continuity in time, we mainly use the result of \cite[Theorem
3]{ish}) with the fact that $u_{t}=-F(x,t,u_{x})$, and the boundedness of 
the Hamiltonian.$\hfill{\Box}$

%%%%%%%%%%%%%%%%%%%%%%%%%%%%%%%%%%%%%%%%%%%%%%%%%%%%%%%%%%%%%%%%%%%%%%%
\begin{rem}\label{uniform_bound}
It is worth mentioning
that the space Lipschitz constant of the function $u$ depends on
$C$, where $C$ appears in $\bold{(F1)}$ for $p=q$, and on
the Lipschitz constant $\gamma$ of the function $u_0$. While the
time Lipschitz constant depends on the bound of the Hamiltonian.
\end{rem}

%%%%%%%%%%%%%%%%%%%%%%%%%%%%%%%%%%%%%%%%%%%%%%%%%%%%%%%%%%%%%%%%%%%%
\subsection{Entropy solution: definition and properties}

%%%%%%%%%%%%%%%%%%%%%%%%%%%%%%%%%%%%%%%%%%%%%%%%%%%%%%%%%%%%%%%%

\begin{definition}\label{def_ent}{\bf(Entropy sub-/super-solution)}\\
Let $F(x,t,v)=g(x,t)f(v)$ with $g,\,g_{x}\in L_{loc}^{\infty}({Q}_T;\R)$
and $f\in
C^{1}(\R;\R)$. A function $v\in L^{\infty}(Q_T;\R)$ is 
an entropy sub-solution of (\ref{CL}) with bounded initial data $v^{0}\in
L^{\infty}(\R)$ if it satisfies:
\begin{equation}\label{subeinq}
\begin{aligned}
\int_{Q_T}\Big[\,\eta_{i}(v(x,t))\phi_{t}(x,t)+\Phi(v(x,t))g(x,t)\phi_{x}(x,t)+\\
h(v(x,t))g_{x}(x,t)\phi(x,t)\,\Big]dxdt+\int_{\R}\eta_{i}(v^{0}(x))\phi(x,0)dx\geq
0,
\end{aligned}
\end{equation}
$\forall\phi\in C^{1}_{0}(\R\times[0,T);\R_{+})$, for any non-decreasing
convex function $\eta_{i}\in C^{1}(\R;\R)$, $\Phi\in C^{1}(\R;\R)$ such
that:
\begin{equation}\label{ent_flow}
\Phi^{'}=f^{'}\eta^{'}_{i}, \quad\mbox{and}\quad h=\Phi-f\eta^{'}_{i}.
\end{equation}
An entropy super-solution of (\ref{CL}) is defined by replacing in
(\ref{subeinq}) $\eta_{i}$ with $\eta_{d}$; a non-increasing convex
function. An entropy solution is defined as being both entropy sub-
and super-solution. In other words, it verifies (\ref{subeinq}) for
any convex function $\eta\in C^{1}(\R;\R)$.
\end{definition}

%%%%%%%%%%%%%%%%%%%%%%%%%%%%%%%%%%%%%%%%%%%%%%%%%%%%%%%%%%%%%%%%%%%%%%

A well know characterization of the entropy solution is that:
\begin{proposition} A function $v\in L^{\infty}(Q_T)$ is an entropy 
sub-solution of (\ref{CL})
if and only if $\forall k\in \R$, $\phi\in C_{0}^{1}(\R\times[0,T);\R_{+})$,
one has:
\begin{equation}
\int_{Q_T}\Big[\,(v(x,t)-k)^{+}\phi_{t}(x,t)+\mbox{sgn}^{+}(v(x,t)-k)(f(v(x,t))-f(k))g(x,t)\phi_{x}(x,t)-
\nonumber
\end{equation}
\begin{equation}\label{withsgn}
\mbox{sgn}^{+}(v(x,t)-k)f(k)g_{x}(x,t)\phi(x,t)\,\Big]dxdt+\int_{R}(v^{0}(x)-k)^{+}\phi(x,0)dx\geq
0,
\end{equation}
Where $a^{\pm}= \frac{1}{2}(|a|\pm a)$ and
$\mbox{sgn}^{\pm}(x)=\frac{1}{2}(\mbox{sgn}(x)\pm 1)$. An entropy 
super-solution of
(\ref{CL}) is defined replacing in (\ref{withsgn})
$(\cdot)^+$, $\mbox{sgn}^+$ by $(\cdot)^-$, $\mbox{sgn}^-$.
\end{proposition}
This characterization can be deduced from (\ref{subeinq}), by using
regularizations of the function $(\cdot-k)^{+}$. Also
(\ref{subeinq}) may be obtained from (\ref{withsgn}) by
approximating any non-decreasing convex function $\eta_{i}\in
C^{1}(\R;\R)$ by a sequence of functions of the form:
$\eta^{(n)}_{i}(\cdot)=\sum^{n}_{1}\beta^{(n)}_{i}(\cdot-k^{(n)}_{i})^{+}$,
with $\beta^{(n)}_{i}\geq 0$.

Entropy solution was first introduced by Kru$\check{z}$kov \cite{kru1} as
the only  physically admissible solution among all weak (distributional)
solutions to scalar conservation laws. These weak solutions lack the fact of
being unique for it is easy to construct multiple weak solutions to Cauchy
problems (\ref{CL}), see \cite{lax}.

Our next definition concerns classical sub-/super-solution to scalar
conservation laws. This kind of solutions are shown to be entropy
solutions, for the details see lemma \ref{subcsube}.
\begin{definition}{\bf(Classical solution to scalar conservation laws)}\\
Let $F(x,t,v)=g(x,t)f(v)$ with  $g,\,g_{x}\in 
L_{loc}^{\infty}(Q_T;\R)$ and $f\in
C^{1}(\R;\R)$. A function $v\in W^{1,\infty}(Q_T)$ is said to be a classical 
sub-solution of
(\ref{CL}) with $v^{0}(x)=v(x,0)$ if it satisfies
\begin{equation}\label{alijinan}
v_{t}(x,t)+(F(x,t,v(x,t)))_{x}\leq 0 \quad\mbox{a.e. in}\quad Q_T.
\end{equation}
Classical super-solutions are defined by replacing ``$\leq$'' with ``$\geq$''
in (\ref{alijinan}), and classical solutions are defined to be both
classical sub- and super-solutions.
\end{definition}

We move now to some results on entropy solutions depicted from
\cite{kru1}.
\begin{theorem}\label{kru}{\bf(Kru$\check{z}$kov's Existence Theorem)}\\
Let $F$, $v^0$ be given by Definition \ref{def_ent}, and the following 
conditions hold:\\\\
$\bold{(E0)}$ $f\in C^{1}_{b}(\R),$\\\\
$\bold{(E1)}$ $g,\,g_{x}\in C_{b}(\bar{Q}_T),$\\\\
$\bold{(E2)}$ $g_{xx}\in C(\bar{Q}_T)$,\\\\
then there exists an entropy solution $v\in
L^{\infty}(Q_T)$ of (\ref{CL}).
\end{theorem}
In fact, Kru$\check{z}$kov's conditions for existence were given for a
general flux function \cite[Section 4]{kru1}. However, in
Subsection 5.4 of the same paper, a weak version of these conditions, that 
can be easily checked in the case  $F(x,t,v)=g(x,t)f(v)$ and
$\bold{(E0)}$-$\bold{(E1)}$-$\bold{(E2)}$, is
presented. Furthermore, uniqueness follows from the following comparison 
principle.

\begin{theorem}\label{entcomp}{\bf(Comparison Principle)}\\
Let $F$ be given by Definition \ref{def_ent} with $f$ satisfying {\bf
  (E0)}, and $g$ satisfies,\\\\
$\bold{(E3)}$ $g\in W^{1,\infty}(\bar{Q}_T).$\\\\
Let $u(x,t)$, $v(x,t) \in L^{\infty}({Q}_T)$ be two entropy sub-/super-solutions of (\ref{CL}) with initial
data $u^{0}$, $v^{0} \in L^{\infty}(\R)$. Suppose that,
$$u^{0}(x)\leq v^{0}(x)\;\;\mbox{ a.e. in } \;\;\R,$$
then
$$ u(x,t)\leq v(x,t) \;\;\mbox{ a.e. in } \;\;\bar{Q}_T.$$
\end{theorem}
\noindent {\bf Proof.} See Section \ref{sec6}, Appendix. $\hfill\Box$\\

It is worth noticing that in \cite{kru1}, the proof of the existence of
entropy solutions of (\ref{CL})  is made through a parabolic
regularization of (\ref{CL}) and passing to the limit, with respect to
the $L^1$ convergence on compacts, in a convenient space.

At this stage, we are ready to present a relation that sometimes
hold between scalar conservation laws and Hamilton-Jacobi equations in
one-dimensional space.

%%%%%%%%%%%%%%%%%%%%%%%RELATION E-V%%%%%%%%%%%%%%%%%%%%%%%%%%%%%%%%
\subsection{Entropy-Viscosity relation}

Formally, by differentiating (\ref{HJ1}) with respect to $x$ and
defining $v=u_x$, we see that (\ref{HJ1}) is equivalent to the
scalar conservation law (\ref{CL}) with $v^{0}=u^{0}_{x}$ and the
same $F$. This equivalence of the two problems has been exploited in
order to translate some numerical methods for hyperbolic
conservation laws to methods for Hamilton-Jacobi equations.
Moreover, several proofs were given in the one dimensional case. The
usual proof of this relation depends strongly on the known results
about existence and uniqueness of the solutions of the two problems
together with the convergence of the viscosity method (see
\cite{corfalnat,kru2,lions1}). Another proof of this relation could be found in
\cite{cas1} via the definion of viscosity/entropy inequalities, while a
direct proof could also be found in 
\cite{karres1} using the front tracking method. The case of a
Hamiltonian of the form (\ref{shape_hamil}) is also treated even
when $g(x,t)$ is allowed to be discontinuous in the $(x,t)$ plane
along a finite number of (possibly intersected) curves, see
\cite{dostrov}.

In our work, the above stated relation will be successfully used to get some
gradient estimates of $\kappa$. Although several approaches were given
to establish this connection, we will present for the reader's
convenience, a proof similar to that given in \cite[Theorem
2.2]{corfalnat}. For every Hamiltionian/flux function $F=gf$ and every $u^{0}\in
Lip(\R)$, let
$$\mathcal{E}\mathcal{V}=\{{\bf (V0)}, {\bf (V1)}, {\bf
    (V2)}, {\bf (E0)}, {\bf (E1)}, {\bf (E2)}, {\bf (E3)}\},$$
in other words,
\begin{equation}
\mathcal{E}\mathcal{V}=\left|
\begin{aligned}
&\mbox{The set of all conditions on $f$ and $g$ ensuring the}\\
&\mbox{existence and uniqueness of a Lipschitz continuous
  viscosity}\\
&\mbox{solution $u\in Lip(\bar{Q}_T)$ of
(\ref{HJ1}), and of an entropy }\\
&\mbox{solution $v\in
L^{\infty}(Q_T)$ of (\ref{CL}), with $v^{0}=u^{0}_{x}\in
L^{\infty}(\R)$.}
\end{aligned}
\right.\nonumber
\end{equation}
\begin{theorem} \label{relation}{\bf(A link between viscosity and entropy 
solutions)}\\
Let $F=gf$ with $g\in C^{2}(\bar{Q}_T)$, $u^{0}\in Lip(\R)$ and
$\mathcal{E}\mathcal{V}$ satisfied. Then,
$$v=u_{x}\quad \mbox{a.e. in}\quad Q_T.$$
\end{theorem}
\noindent{\bf Sketch of the proof.} Let $\varepsilon>0$
and $\delta>0$. We start the prove by making a parabolic
regularization of equation (\ref{HJ1}) and a smooth regularization
of $u_{0}$ and we solve the following parabolic equation:
\begin{equation}\label{epsdel}
\left\{
\begin{aligned}
&u^{\varepsilon,\delta}_{t}+F(x,t,u^{\varepsilon,\delta}_{x})=\epsilon
u^{\varepsilon,\delta}_{xx}\quad&\mbox{in}&\quad \mathbb{R}\times(0,T),\\
&u^{\varepsilon,\delta}(x,0)=u^{0,\delta}(x)\quad&\mbox{in}&\quad 
\mathbb{R}.
\end{aligned}
\right.
\end{equation}
For the sake of simplicity, we will denote $u^{\varepsilon,\delta}$ by $w$ 
and
$u^{0,\delta}$ by $w^{0}$. Note that the first equation of
(\ref{epsdel}) can be viewed as the heat equation with a source term
$F$. Thus, we have:
\begin{equation}\label{heatequation}
\left\{
\begin{aligned}
&w_{t}-\varepsilon w_{xx}=F[w](x,t)\quad&\mbox{in}&\quad Q_T,\\
&w(x,0)=w^{0}\quad&\mbox{in}&\quad \R,
\end{aligned}
\right.
\end{equation}
with $F[w](x,t)=F(x,t,w_x(x,t))$. From the classical theory of heat
equations, since $F[w]\in L_{loc}^{p}(Q_T)$ and $w^{0}\in
W_{loc}^{1,p}(\R)$, there exists a unique solution $w$ of
(\ref{heatequation}) such that
$$  w\in W_{p}^{2,1}(\Omega)\quad\forall\Omega\subset\subset Q_T 
\;\;\mbox{and}\; \;1<p<\infty.$$
Here the space $W_{p}^{2,1}(\Omega)$, $p\geq 1$ is the
Banach space consisting of all functions $w\in L^{p}(\Omega)$ having
generalized derivatives of the form
$w_t$ and $w_{xx}$ in $L^{p}(\Omega)$. For more details, see \cite[Theorem
9.1]{LSU}. We also notice that the space $W_{p}^{2,1}(\Omega)$ is
continuously injected in the Hölder space
$C^{\alpha,\alpha/2}(\Omega)$ for $\alpha=2-\frac{3}{p}$ and
$p>\frac{3}{2}$, see \cite{LSU}. We use now a
bootstrap argument to increase the regularity of $w$, taking in
each stage, the new regularity of $F[w]$ and the regularity of
$w^{0}$. Finally, we get that $w\in C^{3,1}(\R\times[0,T))$ (three
times continuously differentiable in space and one time
continuously differentiable in time). From the maximum principle and the 
$L^p$-estimates of
the heat equation, see \cite{LSU,brez1}, it follows the uniform
bound of $u^{\varepsilon,\delta}$ in $W^{1,p}_{loc}(Q_T)$, for
$p>2$. Therefore, we get as $\delta\rightarrow 0$ and
$\varepsilon\rightarrow 0$ that:
$$u^{\varepsilon,\delta}\rightarrow u \quad\mbox{in}\quad 
C(\R\times[0,T)),$$
with $u(x,0)=u^{0}$. We now make use of the stability theorem, 
\cite[Théorème 2.3]{barbook}, twice
on the equation (\ref{epsdel}) to get that the limit $u$ is the unique
viscosity solution of (\ref{HJ1}). Hence, we have for any 
$\phi\in C^{\infty}_{0}(Q_T)$
\begin{eqnarray*}
\lim_{\varepsilon\rightarrow
  0,\,\delta\rightarrow 
0}\int_{0}^{T}\int_{\mathbb{R}}u_{x}^{\varepsilon,\delta}\phi\, dx\,dt&=&
-\lim_{\varepsilon\rightarrow
  0,\,\delta\rightarrow 
0}\int_{0}^{T}\int_{\mathbb{R}}u^{\varepsilon,\delta}\phi_{x}\,dx\,dt\\
&=&-\int_{0}^{T}\int_{\mathbb{R}}u\phi_{x}\,dx\,dt=\int_{0}^{T}\int_{\mathbb{R}}u_{x}\phi\,
dx\,dt.
\end{eqnarray*}
The appearance of $u_x$ follows since $u\in Lip(\bar{Q}_T)$. Moreover,
as a regular solution, the function 
$v^{\varepsilon,\delta}=u^{\varepsilon,\delta}_x$ solves
the derived problem
\begin{equation}
\left\{
\begin{aligned}
&v^{\varepsilon,\delta}_{t}+(F(x,t,v^{\varepsilon,\delta}))_{x}=\epsilon
v^{\varepsilon,\delta}_{xx}\quad&\mbox{in}&\quad \mathbb{R}\times(0,T),\\
&v^{\varepsilon,\delta}(x,0)=u^{0,\delta}_{x}(x)\quad&\mbox{in}&\quad 
\mathbb{R},
\end{aligned}
\right.
\end{equation}
and, according to \cite[Theorem 4]{kru1}, the sequence 
$v^{\varepsilon,\delta}$ converge in
$L^{1}_{loc}(\bar{Q}_T)$, as $\varepsilon\rightarrow 0$ and
$\delta\rightarrow 0$, to the entropy solution
$v$ of (\ref{CL}). Then, for any $\phi\in C^{\infty}_{0}(Q_T)$,
$$ \lim_{\varepsilon\rightarrow
  0\,\delta\rightarrow 
0}\int_{0}^{T}\int_{\mathbb{R}}v^{\varepsilon,\delta}\phi\, 
dx\,dt=\int_{0}^{T}\int_{\mathbb{R}}v\phi\, dx\,dt.$$
Consequently,
$$\int_{0}^{T}\int_{\mathbb{R}}u_{x}\phi\, 
dx\,dt=\int_{0}^{T}\int_{\mathbb{R}}v\phi\,
dx\,dt,$$
and $u_x=v$ a.e. in $Q_T$. $\hfill\Box$

\begin{rem}
The converse of the previous theorem holds under certain assumptions
(see \cite{karres1, cocris1}).
\end{rem}
\begin{rem}
In the multidimensional case this one-to-one correspondence no
longer exists, instead the gradient $v=\nabla u$ satisfies formally
a non-strict hyperbolic system of conservation laws (see
\cite{lions1, kru2}).
\end{rem}
\end{section}
Throughout Sections 3 and 4, $\rho$ will always be the solution
of the heat equation (\ref{mainheat}). The properties of the solution of the
heat equation with such a regular initial data will be frequently
used, we refer the reader to \cite{brez1,evans} for details.
%*****************************************************************************
%%%%%%%%%%%%%%%%%%%%%%%%%%%%%%%APPROXIMATED PROBLEM%%%%%%%%%%%%%%%%%%%%%%%%%%%

\begin{section}{The approximate problem}\label{sec3}
 In this section, we approximate (\ref{mainhj}) and we pose a more restrictive
condition (see condition (\ref{res_cond_ini_vis})) on the gradient of the initial data than
of the physicaly relevent one (\ref{condition}). We prove a result of existence and uniqueness of
this approximate problem, namely Theorem \ref{main1}, and
the reader will notice at the end of this section that this restrictive condition is satisfied for
all time, and this what cancels the approximation in the structure of (\ref{mainhj}) and
returns it to its original one. Finaly we present the proof of
Theorem \ref{main2}.

For every $a>0$, we build
up an approximation function $f_a\in C^{1}_{b}(\R)$ of the function
$\frac{1}{x}$ defined by:
\begin{equation}\label{functionf}
f_{a}(x) =
\left\{
\begin{aligned}
&\frac{1}{x}\quad \quad\quad\mbox{if} \quad\quad\quad x\geq a,\\
&\frac{2a-x}{a^{2}+a^{2}(x-a)^{2}}\quad \mbox{otherwise.}
\end{aligned}
\right.
\end{equation}
\begin{proposition}\label{approblem}
For any $a>0$, let $f_a$ be defined by (\ref{functionf}) and $H\in
C^{1}(\R)$ be a scalar-valued function. If
\begin{equation}\label{fdef}
F_{a}(x,t,u)=-H(\rho_{x}(x,t))\rho_{xx}(x,t)f_a(u)
\end{equation}
and $\kappa^{0}\in Lip(\mathbb{R})$, then the Hamilton-Jacobi equation
\begin{equation}\label{hjapp}
\left\{
\begin{aligned}
&\kappa_{t}+F_{a}(x,t,\kappa_x)=0\quad&\mbox{in}&\quad Q_T,\\
&\kappa(x,0)=\kappa^{0}(x)\quad&\mbox{in}&\quad \mathbb{R},
\end{aligned}
\right.
\end{equation}
has a unique viscosity solution $\kappa\in Lip(\bar{Q}_T)$.
\end{proposition}
\noindent{\bf Proof.} The proof is easily concluded from Theorems
\ref{comphj}, \ref{exishj} and Proposition \ref{lipreghj}, after
checking that the conditions
$\bold{(V0)}$-$\bold{(V1)}$-$\bold{(V2)}$ are satisfied with
\begin{equation}\label{fun_g}
g(x,t)=-H(\rho_{x}(x,t))\rho_{xx}(x,t).
\end{equation}
The condition $\bold{(V0)}$ is
trivial, while for $\bold{(V1)}$, we just use the fact that $H$ is bounded
on compacts and the fact that $|\rho_{x}(x,t)|\leq
||\rho^{0}_{x}||_{L^{\infty}(\R)}$ in $\bar{Q}_T$. For the condition
$\bold{(V2)}$, the regularity of $\rho$ and $H$ permits to compute the
spatial derivative of $g$ in $\bar{Q}_T$, thus we have:
$$g_{x}=-(H^{'}(\rho_{x})\rho_{xx}^{2}+H(\rho_{x})\rho_{xxx}).$$
The uniform bound of the spatial derivatives, up to the third order, of the 
solution of the heat
equation, and the boundedness of $H^{'}$ on compacts gives immediately
$\bold{(V2)}$. $\hfill\Box$\\

In the following proposition, we show a lower-bound estimate for the gradient of 
$\kappa$ obtained in Proposition \ref{approblem}. It is worth mentioning
that a result of lower-bound gradient estimates for first-order
Hamilton-Jacobi equations could be found in \cite[Theorem
4.2]{ley}. However, this result holds for Hamiltonians $F(x,t,u)$ that are convex
in the $u$-variable, using only the viscosity theory techniques. This is
not the case here, and in order to obtain our lower-bound estimates, we
need to use the viscosity/entropy theory techniques. In particular, we
have the following:
\begin{proposition}\label{addpro}
Let $G\in C^{3}(\mathbb{R};\R)$ satisfying the following conditions:
\begin{itemize}
\item[(G1)] $G(x)\geq G(0)>0$,
\item[(G2)] $G^{''}\geq 0.$
\end{itemize}
Moreover, let $$H=GG^{'}\quad \mbox{and} \quad 0<a\leq G(0).$$
If $\kappa^{0}$ satisfies:
$$\kappa^{0}_{x}(x)\geq G(\rho^{0}_{x}(x)),\quad \mbox{a.e. in }\quad\R,$$
then the solution $\kappa$ obtained from Proposition \ref{approblem} 
satisfies:
 \begin{equation}\label{maininq}
\kappa_{x}(x,t)\geq G(\rho_{x}(x,t))\quad \mbox{ a.e.}\;\mbox{ in
}\quad \bar{Q}_T.
\end{equation}
\end{proposition}
In order to prove Proposition \ref{addpro}, we first show
that $G(\rho_{x})$ is an entropy sub-solution of
\begin{equation}\label{dereqn}
\left\{
\begin{aligned}
&\omega_{t}+(F(x,t,\omega))_{x}=0 \quad&\mbox{in}&\quad Q_T,\\
&\omega(x,0)=\omega^{0}(x)\quad&\mbox{in}&\quad\mathbb{R},
\end{aligned}
\right.
\end{equation}
with $w^{0}=G(\rho^{0}_{x})$ and $F$ is the same as in (\ref{fdef}). Before
going further, we will pause to prove a lemma which makes it easier to reach 
our goal.
\begin{lemma}\label{subcsube}{\bf (Classical sub-solutions are entropy 
sub-solutions)}\\
Let $v\in W^{1,\infty}(Q_T)$ be a classical sub-solution of (\ref{CL})
with $v^{0}(x)=v(x,0)$, then $v$ is an entropy sub-solution.
\end{lemma}
\noindent{\bf Proof.} Let $\eta_{i}$, $\Phi$, $h$ and
$\phi$ be given by Definition \ref{def_ent}. Multiplying inequality
(\ref{alijinan}) by $\eta^{'}_{i}(v)\phi$ does not change its sign.
Hence, after developing, we have:
\begin{equation}\label{clas1}
\eta^{'}_{i}(v)v_{t}\phi+\eta^{'}_{i}(v)g_{x}f(v)\phi+\eta^{'}_{i}(v)gf^{'}(v)v_{x}\phi\leq
0, \quad\mbox{a.e. in }\;Q_T,
\end{equation}
and since $v$ is Lipschitz continuous, we use the chain-rule formula together with
(\ref{ent_flow}) to rewrite (\ref{clas1}) as:
\begin{equation}\label{clas2}
(\eta_{i}(v))_{t}\,\phi+g_{x}f(v)\eta_{i}^{'}(v)\phi+g(\Phi(v))_{x}\,\phi\leq
0,  \quad\mbox{a.e. in }\;Q_T.
\end{equation}
Upon integrating (\ref{clas2}) over $Q_T$ and transferring derivatives
with respect to $t$ and $x$ to the test function, we obtain:
\begin{equation}\label{clas3}
\hspace{-3.2cm}\int_{Q_T}\Big[\,\eta_{i}(v(x,t))\phi_{t}(x,t)+\Phi(v(x,t))g(x,t)\phi_{x}(x,t)+
\nonumber
\end{equation}
\begin{equation}\label{clas4}
\hspace{1.5cm}h(v(x,t))g_{x}(x,t)\phi(x,t)\,\Big]dxdt+\int_{\R}\eta_{i}(v^{0}(x))\phi(x,0)dx\geq 
0,
\end{equation}
which ends the proof.  $\hfill\Box$\\

\noindent Following same arguments, classical super-solutions are shown to 
entropy
super-solutions. We return now to the function $G(\rho_{x})$ and we are ready 
to show that
it is indeed an entropy sub-solution of (\ref{dereqn}). In particular, we 
have the following:
\begin{lemma}\label{strsubsol}
The function $G(\rho_{x})$ defined on $Q_T$ is a classical 
sub-solution of (\ref{dereqn}) with initial data
$G(\rho^{0}_{x})$, hence an entropy sub-solution.
\end{lemma}
\noindent{\bf Proof of Lemma \ref{strsubsol}.} First, it is easily seen that $G(\rho_{x})\in
W^{1,\infty}(Q_T)$. Define the scalar valued quantity $B$ on $Q_T$
by:
$$B(x,t)=\partial_{t}(G(\rho_{x}(x,t)))+\partial_{x}(F(x,t,G(\rho_{x}(x,t)))).$$
Since $0<a\leq G(0)$, we use (G1) to get
$f_{a}(G(\rho_{x}))=1/G(\rho_{x})$ and we observe that,
\begin{eqnarray*}
B&=&G^{'}(\rho_{x})\rho_{xt}-\partial_{x}\left(\frac{H(\rho_{x})\rho_{xx}}{G(\rho_{x})}\right)\\
&=&G^{'}(\rho_{x})\rho_{xxx}-\left(\frac{G(\rho_{x})[H^{'}(\rho_{x})\rho_{xx}^{2}+H(\rho_{x})\rho_{xxx}]-(G^{'}(\rho_{x})
\rho_{xx}^{2}H(\rho_{x}))}{G^{2}(\rho_{x})}\right)\\
&=&
\frac{G(\rho_{x})\rho_{xxx}
  (G(\rho_{x})G^{'}(\rho_{x})-H(\rho_{x}))-\rho_{xx}^{2}(H^{'}(\rho_{x})G(\rho_{x})-H(\rho
_{x})G^{'}(\rho_{x}))}{G^{2}(\rho_{x})}\\
&=&-\rho_{xx}^{2}G^{''}(\rho_{x}).
\end{eqnarray*}
The condition (G2) gives immediately that $B\leq 0$. This proves
that $G(\rho_{x})$ is a classical sub-solution of equation
(\ref{dereqn}) and hence an entropy sub-solution.
$\hfill\Box$\\

\noindent{\bf Proof of Proposition \ref{addpro}.} From the
definition of $H$ and the properties of $\rho$, it is easy to check that 
$g\in C^{2}(\bar{Q}_T)$ and
that $\mathcal{E}\mathcal{V}$ is fully satisfied. Hence, we are in the
framework of Theorem~\ref{relation} with $u^{0}=\kappa^{0}$. This theorem 
gives that $\kappa_{x}$ is the
unique entropy solution of (\ref{dereqn}) with $w^{0}=\kappa^{0}_{x}$.  
Moreover,  by the previous
lemma, $G(\rho_{x})$ is an entropy sub-solution of (\ref{dereqn}). Since 
$$\kappa^{0}_{x}\geq
G(\rho^{0}_{x}),\quad \mbox{a.e. in}\quad\R,$$ we can apply the Comparison
Theorem \ref{entcomp} to get the desired result.$\hfill\Box$\\

%%%%%%%%%%%%%%%%%%%%%%%%%%%%%%%%%%%%%%%%%%%%%%%%%%%%%%%%%%%%
\noindent It is worth notable here that we do not know how to obtain the
lower-bound on the spatial gradient $\kappa_{x}$ using the viscosity framework
directly. However, for the case of the upper-bound, we can do so
(see Remark \ref{remark_v}). At this stage, fix some $\epsilon>0$, and let
$$G^{\epsilon}(x)=\sqrt{x^{2}+\epsilon^{2}} \quad \mbox{and} \quad
a=G^{\eps}(0)=\eps.$$ It is clear that $G^{\epsilon}(x)$
satisfies the conditions (G1)-(G2) with
$$H^{\epsilon}(x)=x,$$ 
 and the Hamiltonian $F$ from
(\ref{fdef}) takes now the following shape:
\begin{equation}\label{newf_def}
F_{\eps}(x,t,u)=-\rho_{x}(x,t)\rho_{xx}(x,t)f_{\eps}(u).
\end{equation}
Moreover, we have the
following corollary which is is an immediate consequence of Propositions
\ref{approblem} and \ref{addpro}.
\begin{corollary}\label{realapp}
There exists a unique viscosity solution $\kappa\in Lip(\bar{Q}_T)$ of
\begin{equation}\label{realappeqn}
\left\{
\begin{aligned}
&\kappa_{t}+F_{\eps}(x,t,\kappa_{x})=0\quad&\mbox{in}&\quad
Q_T,\\
&\kappa(x,0)=\kappa^{0}\in Lip(\mathbb{R})\quad&\mbox{in}&\quad
\mathbb{R},
\end{aligned}
\right.
\end{equation}
with $\kappa^{0}_{x}$ satisfies:
\begin{equation}\label{res_condition}
\kappa_{x}^{0}\geq
\sqrt{(\rho_{x}^{0})^{2}+\epsilon^{2}}\quad\mbox{
a.e. in}\quad \mathbb{R}.
\end{equation}
Moreover, this solution $\kappa$ satisfies:
\begin{equation}\label{lolo} \kappa_{x}\geq
\sqrt{\rho^{2}_{x}+\epsilon^{2}}\quad\mbox{ a.e.  in }
\quad\bar{Q}_T.
\end{equation}
\end{corollary}
The following lamma will be used in the proof of Theorem \ref{main1}.
\begin{lemma}\label{123}
Let $\bar{c}$ be an arbitrary real constant and take $\psi\in
Lip(\R;\R)$ satisfying:
$$\psi_{x}\geq \bar{c}\quad\mbox{a.e. in}\quad \R.$$
If $\zeta\in C^{1}(\R;\R)$ is such that $\psi-\zeta$ has a local
maximum or local minimum at some point $x_{0}\in\R$, then
$$\zeta_{x}(x_{0})\geq \bar{c}.$$
\end{lemma}
\noindent {\bf Proof.} Suppose that $\psi-\zeta$ has a local minimum
at the point $x_{0}$; this ensures the existence of a certain $r>0$
such that
$$(\psi-\zeta)(x)\geq (\psi-\zeta)(x_{0})\quad \forall x; \;|x-x_{0}|<r.$$
We argue by contradiction. Assuming $\zeta_{x}(x_{0})<\bar{c}$ leads, from the continuity of
$\zeta_{x}$, to the existence of $r^{'}\in (0,r)$ such that
\begin{equation}\label{tayyar}
\zeta_{x}(x)<\bar{c}\quad \forall x;\;|x-x_{0}|<r^{'}.
\end{equation}
Let $y_{0}$ be a point such that $|y_{0}-x_{0}|<r^{'}$ and
$y_{0}<x_{0}$. Reexpressing (\ref{tayyar}), we get
$$(\zeta-\bar{c}x)_{x}(x)<0 \quad \forall x\in(y_{0},x_{0}),$$ and hence
$$\int_{y_{0}}^{x_{0}}[(\psi-\bar{c}x)_{x}(x)-(\zeta-\bar{c}x)_{x}(x)]dx>0,$$
which implies that
$$(\psi-\zeta)(x_{0})>(\psi-\zeta)(y_{0}),$$
and hence a contradiction. We remark that the case of a local
maximum can be treated in a
similar way.$\hfill\Box$\\

Now, we are ready to present the proofs of the first two theorems
announced in Section \ref{sec1}.\\\\
 \noindent{\bf Proof of Theorem \ref{main1}.} Let
$\kappa\in Lip(\bar{Q}_{T})$ be the solution of (\ref{realappeqn})
obtained in Corollary \ref{realapp}. Let us show that it is the unique
viscosity solution of (\ref{mainhj}) among those verifying (\ref{lolo}). To do this, we
consider a test function $\phi\in C^{1}(Q_T)$ 
such that $\kappa-\phi$ has a local minimum at some point
$(x_{0},t_{0})\in Q_T$. Proposition \ref{lipreghj}, together with
inequality (\ref{lolo}) gives that
$$\kappa(.,t_{0})\in Lip(\R)\quad\mbox{and}\quad\kappa_{x}(.,t_{0})\geq\eps
\;\;\mbox{a.e. in } \;\R.$$
We make use of Lemma \ref{123} with $\psi(.)=\kappa(.,t_{0})$ and
$\zeta(.)=\phi(.,t_{0})$ to get
\begin{equation}\label{456}
\phi_{x}(x_{0},t_{0})\geq \eps.
\end{equation}
Since $\kappa$ is a viscosity super-solution of
$$\kappa_{t}-f_{\eps}(\kappa_{x})\rho_{x}\rho_{xx}=0\quad\mbox{in}\quad
Q_T,$$
we have
$$\phi_{t}(x_{0},t_{0})-f_{\eps}(\phi_{x}(x_{0},t_{0}))\rho_{x}(x_{0},t_{0})\rho_{xx}(x_{0},t_{0})\geq
0.$$
However, from (\ref{456}), we get
$$\phi_{t}(x_{0},t_{0})\phi_{x}(x_{0},t_{0})-\rho_{x}(x_{0},t_{0})\rho_{xx}(x_{0},t_{0})\geq
0,$$
and hence $\kappa$ is a viscosity super-solution of
$$\kappa_{t}\kappa_{x}=\rho_{x}\rho_{xx} \quad\mbox{in}\quad
Q_T.$$ In the same way, we can show that $\kappa$ is a viscosity
sub-solution of the above equation and hence a viscosity solution.
The uniqueness of this solution comes from the uniqueness of the
viscosity solution of (\ref{realappeqn}) by reversing the above
reasoning. $\hfill\Box$
\begin{rem}\label{notdepep}
Notice that the first equation of (\ref{mainhj}) can be viewed as a
Hamilton-Jacobi equation of the type
$$F(X,\nabla\kappa)=0\quad \mbox{in} \quad Q_T,$$ where
$F:Q_{T}\times\mathbb{R}^{2}\mapsto \mathbb{R}$
defined by:
$$F(X,p)=p_{1}p_{2}-\rho_{x}(X)\rho_{xx}(X),$$ with $X=(x,t)$ and
$p=(p_{1},p_{2})$.
\end{rem}
\noindent {\bf Proof of Theorem \ref{main2}.} Let $\theta=\kappa_x$. By
Theorem \ref{relation}, $\theta$ is the unique entropy solution of
\begin{eqnarray*}
\left\{
\begin{aligned}
&\theta_{t}=(\rho_{x}\rho_{xx}f_{\eps}(\theta))_{x}\quad&\mbox{
  in}&\quad Q_T,\\
&\theta(x,0)=\theta^{0}(x)\quad&\mbox{in}&\quad\mathbb{R},
\end{aligned}
\right.
\end{eqnarray*}
with
$$
\theta^{0}(x)=\kappa_{x}^{0}(x)\geq\sqrt{(\rho_{x}^{0})^{2}+\epsilon^{2}},\quad\mbox{a.e.
in}\quad\R.$$
Moreover, from Corollary \ref{realapp}, we have
$$ \theta\geq\sqrt{\rho_{x}^{2}+\epsilon^{2}}\quad\mbox{
  a.e.   in }
\quad\bar{Q}_T,$$
from which we deduce that $f_{\eps}(\theta)=\frac{1}{\theta}$ and hence our theorem holds.$\hfill\Box$
\end{section}
%********************************************************************************
%%%%%%%%%%%%%%%%%%%%%%%%%END OF SECTION%%%%%%%%%%%%%%%%%%%%%%%%%%%%%%%%%%%%%%%%%%
%********************************************************************************
%*********************** Main Theorem *******************************************
\begin{section}{Proof of Theorem \ref{maintheorem}}\label{sec4}
We turn our attention now to Theorem \ref{maintheorem}. Let
$0<\epsilon<1$ be a fixed constant and take
\begin{equation}\label{ja}
\kappa^{0,\epsilon}(x)=\kappa^{0}(x)+\epsilon x.
\end{equation}
It is easy to check that the function $\kappa^{0,\epsilon}$ belongs
to $Lip(\mathbb{R})$, and by condition (\ref{condition}) we get for
a.e. $x\in\R$,
\begin{eqnarray*}
\kappa^{0,\epsilon}_{x}(x)&=&\kappa^{0}_{x}(x)+\epsilon,\\
&\geq&\sqrt{(\rho^{0}_{x}(x))^{2}+\epsilon^{2}}.
\end{eqnarray*}
From  Theorem \ref{main1}, there exists a family of viscosity
solutions $\kappa^{\epsilon}\in Lip(\bar{Q}_{T})$ to the initial
value problem (\ref{mainhj}) that satisfy:
$$\kappa^{\epsilon}_{x}\geq \sqrt{\rho^{2}_{x}+\epsilon^{2}}\quad\mbox{ 
a.e.\;\;\;in }
\quad\bar{Q}_T.$$ We will try to  extract a subsequence of
$\kappa^{\epsilon}$ that converges, in a suitable space, to the desired
solution 

\subsection{Gradient estimates.}

Uniform bounds for the space-time gradients of $\kappa^{\epsilon}$ will play 
an essential role
in the determination of our subsequence.\\\\
{\bf\underline{I. $\eps$-uniform upper-bound for $\kappa^{\epsilon}_{t}$.}}\\
Starting with the time gradient, we have for a.e. $(x,t)\in Q_T$:
\begin{equation}\label{1}
\kappa^{\epsilon}_{t}(x,t)\kappa^{\epsilon}_{x}(x,t)=\rho_{x}(x,t)\rho_{xx}(x,t),
\end{equation}
and
\begin{equation}\label{2}
\kappa^{\epsilon}_{x}(x,t)\geq\sqrt{\rho^{2}_{x}(x,t)+\epsilon^{2}}>0
\quad\mbox{a.e. in} \quad \bar{Q}_T.
\end{equation}
If $\rho_{x}(x,t)=0$ for some Lebesgue point $(x,t)$ of
$\kappa^{\epsilon}_{x}$ and $\kappa^{\epsilon}_{t}$, it
follows from (\ref{1}) and (\ref{2})
that $\kappa^{\epsilon}_{t}(x,t)=0$. Otherwise, and since  by (\ref{2})
$\kappa^{\epsilon}_{x}\geq |\rho_{x}|$, we conclude that:
\begin{equation}\label{grad_es_1}
|\kappa^{\epsilon}_{t}|\leq
||\rho^{0}_{xx}||_{L^{\infty}(\R)}\quad\mbox{a.e. in}\quad Q_T,
\end{equation}
and hence we obtain an $\eps$-uniform bound of
$\kappa^{\epsilon}_{t}$. 

For the space gradient, we argue in a
slightly different way. The key point for obtaining the uniform
bound of $\kappa^{\epsilon}_{t}$ was the minoration of
$\kappa^{\epsilon}_{x}$ by $|\rho_x|$ so, roughly speaking, if we
want to follow the same previous steps using the symmetry of
(\ref{1}) in $\kappa^{\epsilon}_{t}$ and $\kappa^{\epsilon}_{x}$,
one should also have an appropriate minoration of
$|\kappa^{\epsilon}_{t}|$ by a well controlled function
which no longer exists.\\\\
{\bf\underline{II. Formal calculus and best candidate.}}\\
We seek to find the best candidate to be an upper-bound of
$\kappa^{\epsilon}_{x}$. For this reason, we regard formally what is
happening at the maximum of $\kappa^{\epsilon}_{x}$. Dividing both
sides of (\ref{1}) by $\kappa^{\epsilon}_{x}$ and differentiating
with respect to the spatial variable, we get:
\begin{equation}\label{3}
\kappa^{\epsilon}_{xt}=\frac{\rho^{2}_{xx}+\rho_{x}\rho_{xxx}}
{\kappa^{\epsilon}_{x}}-\frac{\kappa^{\epsilon}_{xx}\rho_{x}\rho_{xx}}{(\kappa^{\epsilon}_{x})^{2}}.
\end{equation}
Notice that $\kappa^{\epsilon}_{xx}=0$ at the maximum of
$\kappa^{\epsilon}_{x}$. Multiplying equality (\ref{3}) by
$\kappa^{\epsilon}_{x}$ and integrating between $0$ and $t$, we
obtain:
$$\int^{t}_{0}\frac{d}{d\tau}\left(\frac{1}{2}(\kappa_{x}^{\epsilon})^{2}\right)d\tau=\int^{t}_{0}
(\rho^{2}_{xx}+\rho_{x}\rho_{xxx})d\tau,$$
then
$$(\kappa_{x}^{\epsilon}(x,t))^{2}=(\kappa_{x}^{0,\epsilon}(x))^{2}+2\int^{t}_{0}
(\rho^{2}_{xx}(x,t)+\rho_{x}(x,t)\rho_{xxx}(x,t))d\tau,$$
and hence,
\begin{equation}
|\kappa^{\epsilon}_{x}|\leq \sqrt{2c_{1}t+c_{2}},
\nonumber
\end{equation}
where
$$
c_{1}=||(\rho^{0}_{xx})^{2}||_{L^{\infty}(\R)}+||\rho^{0}_{x}||_{L^{\infty}(\R)}||\rho^{0}_{xxx}||_{L^{\infty}(\R)},$$
and
$$c_{2}=(||\kappa^{0}_{x}||_{L^{\infty}(\R)}+1)^{2}.$$
The reason of taking $c_2$ as above easily follows since
$\kappa^{0,\epsilon}_{x}=\kappa^{0}_{x}+\epsilon$, by taking $\epsilon$
small enough, namely less than $1$.\\\\
{\bf\underline{III. $\eps$-uniform upper-bound for $\kappa^{\epsilon}_{x}$.}}\\
Define the function $S$ by:
$$ S(x,t)=\sqrt{2c_{1}t+c_{2}}.$$
Let us show that $S$ is an entropy super-solution of (\ref{dereqn}) with
$F$ given by (\ref{newf_def}) and 
$w^{0}(x)=S(x,0)$. Indeed, it remark that $S\in W^{1,\infty}(Q_T)$,
and we know that for every $(x,t)\in Q_T$ we have,
$$S(x,t)\geq \sqrt{c_{2}}=||\kappa^{0}_{x}||_{L^{\infty}(\R)}+1\geq\eps,$$
then
\begin{equation}\label{sis}
f_{\eps}(S(x,t))=\frac{1}{S(x,t)} \quad\forall(x,t)\in Q_T.
\end{equation}
The regularity of the function $S$ permits to inject it directly
into the first equation of (\ref{dereqn}). Therefore, using
(\ref{sis}), we have
\begin{eqnarray*}
S_{t}-\left(\frac{\rho_{x}\rho_{xx}}{S}\right)_{x}&=&\frac{c_1}{\sqrt{2c_{1}t+c_{2}}}-\frac{\rho^{2}_{xx}+\rho_{x}\rho_{xxx}}{\sqrt{2c_{1}t+c_{2}}},\\
&=&\frac{c_{1}-(\rho^{2}_{xx}+\rho_{x}\rho_{xxx})}{\sqrt{2c_{1}t+c_{2}}},\\
&\geq&0,
\end{eqnarray*}
which proves, by Lemma \ref{subcsube}, that $S$ is an entropy
super-solution of (\ref{dereqn}). From the discussion of the proof of 
Proposition \ref{addpro},
we know that $\kappa^{\epsilon}_{x}$ is an entropy solution of
(\ref{dereqn}) hence an entropy sub-solution. Since for $\epsilon<1$
and a.e. $x\in\R$, we have,
\begin{eqnarray*}
\kappa^{0,\epsilon}_{x}(x)&=&\kappa^{0}_{x}(x)+\epsilon,\\
&\leq& ||\kappa^{0}_{x}||_{L^{\infty}(\R)}+1,\\
&\leq& \sqrt{c_2}=S(x,0),
\end{eqnarray*}
then we can use the Comparison Theorem \ref{entcomp} of scalar
conservation laws to obtain:
\begin{equation}\label{grad_es_2}
\kappa^{\epsilon}_{x}(x,t)\leq \sqrt{c_{1}t+c_{2}}\leq
\sqrt{c_{1}T+c_{2}}\quad\mbox{a.e. in}\quad \bar{Q}_T,
\end{equation}
and hence we get an $\eps$-uniform bound for $\kappa^{\eps}_{x}$.
\begin{rem}\label{remark_v}
We were able to obtain this $\eps$-uniform upper-bound of
$\kappa^{\eps}_{x}$ by using the viscosity theory techniques. In fact, we claim
that $\zeta^{1,\eps}(x,y,t)=\kappa^{\eps}(x,t)-\kappa^{\eps}(y,t)$ and
 $\zeta^{2}(x,y,t)=(x-y)S(t)$ are two viscosity sub-/super-solutions of
the following Hamilton-Jacobi equation:
\begin{equation}
\frac{\partial w}{\partial t}=F(x,t,w_{x})-F(y,t,-w_{y})
\quad\mbox{in}\quad \mathcal{D}=\{(x,y,t);\;x>y\;\;\mbox{and}\;\;t>0\}
\nonumber
\end{equation}
with initial data
$\zeta^{1,\eps}(x,y,0)=\kappa^{0,\eps}(x)-\kappa^{0,\eps}(y)$ and
$\zeta^{2}(x,y,0)=(x-y)S(0)$ respectively. Here $F$ is given by
(\ref{newf_def}). The claim is easy for $\zeta^{2}$, and we refer to
\cite{CL1} when $\kappa^{\eps}$ is a continuous viscosity solution
of (\ref{realappeqn}). We also notice that:
$\zeta^{1,\eps}(x,y,0)\leq
\zeta^{2}(x,y,0)\;\forall(x,y,0)\in\mathcal{D}$, and
$\zeta^{1,\eps}(x,y,t)=\zeta^{2}(x,y,t)=0$ for $x=y$, $t\geq 0$.
Moreover, since $\zeta^{1,\eps}$ and $\zeta^{2}$ are continuous
functions, we use the comparison principle of viscosity solutions
(see \cite{barbook}) to obtain:
$$\kappa^{\eps}(x,t)-\kappa^{\eps}(y,t)\leq (x-y)S(t)\quad\forall
(x,y,t)\in \bar{\mathcal{D}},$$
hence, the estimate (\ref{grad_es_2}) holds.
\end{rem}
\subsection{Local boundedness in $W^{1,\infty}$.}
We now show that the family $(\kappa^{\epsilon})_{0<\epsilon<1}$ is
locally bounded in $W^{1,\infty}(Q_T)$. Let $K\hspace{-0.1cm}\subset\subset 
Q_T$ be a
compactly contained subset of $Q_T$, and $(x,t)\in K$. Since
$\kappa^{\epsilon}$ is Lipschitz continuous, we can write,
$$ |\kappa^{\epsilon}(x,t)-\kappa^{0,\epsilon}(0)|\leq
C^{\epsilon}_{lip}\,|(x,t)|,$$ where $C^{\epsilon}_{lip}$ is the
Lipschitz constant of $\kappa^{\epsilon}$ which is independent of
$\epsilon$ from the previous estimates, namely (\ref{grad_es_1}) and
(\ref{grad_es_2}). Call this constant $\bar{C}$. From the
definition of $\kappa^{0,\epsilon}(0)$ given by (\ref{ja}), it
follows that,
\begin{eqnarray*}
|\kappa^{\epsilon}(x,t)|&\leq&\bar{C}\,|(x,t)|+|\kappa^{0}(0)|,\\
&\leq&\bar{C}\max_{(y,\tau)\in K}|(y,\tau)|+|\kappa^{0}(0)|,
\end{eqnarray*}
which is finite since $K$ is bounded and hence,
$(\kappa^{\epsilon})_{0<\epsilon<1}$ is uniformly bounded in $C(K)$.
This, together with the uniform gradient estimates, gives the local 
boundedness of $\kappa^{\epsilon}$ in
$W^{1,\infty}(\bar{Q}_T)$.
\subsection{Proof of theorem \ref{maintheorem}}
At this point, we have the necessary tools to give the proof of
Theorem \ref{maintheorem}. We first recall that $\kappa^{\epsilon}$
is a viscosity solution of an equation of the type (\ref{1}), with a
Hamiltonian independent of $\epsilon$ (see Remark \ref{notdepep})
and $\kappa^{0,\epsilon}\rightarrow\kappa^{0}$ locally uniformly in
$\mathbb{R}$. By Ascoli's Theorem, there is a subsequence, called again
$\kappa^{\eps}$, that converges to
$\kappa \in Lip(\bar{Q}_T)$ locally uniformly,
and by the stability theorem (see \cite[Theorem 2.3]{barbook}), $\kappa$ is 
a viscosity solution of the initial value problem
\begin{equation}\label{mainhj1}
\left\{
\begin{aligned}
&\kappa_{t}\kappa_{x}=\rho_{x}\rho_{xx}\hspace{0.64cm} \quad\mbox{ in }\quad 
Q_T,\\
&\kappa(x,0)=\kappa^{0}(x)\hspace{0.31cm}\quad\mbox{ in }\quad \mathbb{R}.
\end{aligned}
\right.
\end{equation}
To end the proof, we still have to show the inequality
$$ \kappa_{x}\geq|\rho_{x}|\quad\mbox{ a.e. in }\quad \bar{Q}_T.$$
Again by Theorem \ref{main1}, our $\kappa^{\epsilon}$ verifies for
a.e. $(x,t)\in \bar{Q}_T$,
\begin{eqnarray*}
\kappa^{\epsilon}_{x}(x,t)&\geq& \sqrt{\rho_{x}^{2}(x,t)+\epsilon^{2}}\\
&>&|\rho_{x}(x,t)|,
\end{eqnarray*}
then for $(y,t)$, $(x,t)\in Q_T$ close enough, with $\rho_{x}$ a
continuous function, the following
inequality hold
$$
\frac{\kappa^{\epsilon}(y,t)-\kappa^{\epsilon}(x,t)}{x-y}>|\rho_{x}(x,t)|.$$
Using the local uniform convergence of $\kappa^{\epsilon}$ to $\kappa$, we
get a similar inequality with $\kappa^{\epsilon}$ replaced with
$\kappa$ and hence
$$\kappa_{x}\geq|\rho_{x}|\quad\mbox{ a.e. in }\quad\bar{Q}_T.$$ 
$\hfill\Box$
\end{section}
%%%%%%%%%%%%%%%%%%%Problem with boundary conditions%%%%%%%%%%%%%%%%%%%%%%%%%%%%%%%%
%%%%%%%%%%%%%%%%%%%%%%%%%%%%%%%%%%%%%%%%%%%%%%%%%%%%%%%%%%%%%%%%%%%%%%%%%%%%%%%%%%%
\begin{section}{Problem with boundary conditions}\label{sec5}
In this part of the paper, we deal with the same problem
structure but with boundary conditions of the Dirichlet type. This
sort of boundary conditions arises naturally in a special model of
dislocation dynamics and will be explained in the following subsection. Our notations are
kept untouched; the terms $\theta^+$, $\theta^-$, $\rho$ and
$\kappa$ still have the same physical meaning, while the domain is
changed into the open and bounded interval
$$I=(0,1),$$
of the real line. Although this problem seems to be an independent one, we 
will
try to benefit the results of the previous sections by considering a
trick of extension and restriction, in order to apply some of the previous
results of the whole space problem.
\subsection{Brief physical motivation}
To illustrate some physical motivations of the boundary value
problem, we consider a constrained channel deforming in simple shear
(see \cite{gcz}). A channel of width $1$ in the $x$-direction and
infinite extension in the $y$-direction is bounded by walls that are
impenetrable for dislocations (see Figure 1). The motion of the
positive and negative dislocations corresponds to the $x$-direction.
\begin{figure}[h!]\label{fig}
\centering\epsfig{file=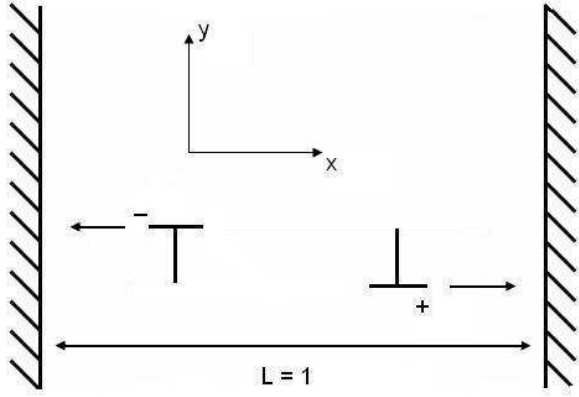, width=100mm}
\caption{Geometry of a constrained channel}
\end{figure}
This is a simplified version of a system studied by Van der Giessen
and coworkers \cite{cvn}, where the simplifications stem from the
fact that:
\begin{itemize}
\item[$\bullet$] only a single slip system is assumed to be active,
such that reactions between dislocations of different type need not
be considered;
\item[$\bullet$] the boundary conditions reduce to "no flux"
conditions for the dislocation fluxes at the boundary walls.
\end{itemize}
The mathematical formulation of this model, as expressed in
\cite{gcz}, is the system (\ref{gcz}) posed on $I\times(0,T)$:
\begin{equation}\label{boundsys}
\left\{
\begin{aligned}
&\partial_{t}\theta^{+}(x,t)-\partial_{x}\left(\theta^{+}(x,t)\left(\frac{\theta^{+}_{x}(x,t)-\theta^{-}_{x}(x,t)}{\theta^{+}(x,t)+\theta^{-}(x,t)}\right)\right)=0,\\
&\partial_{t}\theta^{-}(x,t)+\partial_{x}\left(\theta^{-}(x,t)\left(\frac{\theta^{+}_{x}(x,t)-\theta^{-}_{x}(x,t)}{\theta^{+}(x,t)+\theta^{-}(x,t)}\right)\right)=0.
\end{aligned}
\right.
\end{equation}
To formulate heuristically the boundary conditions at the walls located at
$x=0$ and $x=1$, we note that the dislocation fluxes at the walls
must be zero, which requires that
\begin{equation}\label{flux}
\overbrace{\partial_{x}(\theta^{+}-\theta^{-})}^{\Phi}=0,\qquad\mbox{at}\qquad
x\in \{0,1\}.
\end{equation}
Rewriting system (\ref{boundsys}) in a special integrated form in terms
of $\rho$, $\kappa$ and $\Phi$, we get
\begin{equation}\label{atthebound}
\left\{
\begin{aligned}
&\kappa_{t}=(\rho_{x}/\kappa_{x})\Phi,\\
&\rho_{t}=\Phi.
\end{aligned}
\right.
\end{equation}
Using (\ref{flux}) into the system (\ref{atthebound}), we can
formally deduce that $\rho$ and $\kappa$ are constants along the
boundary walls. Therefore, the remaining of this paper focuses attention on the
study of the following coupled Dirichlet boundary problems:
\begin{equation}\label{heatbound}
\left\{
\begin{aligned}
&\rho_{t}=\rho_{xx}, &\mbox{in}&&I\times(0,\infty),\\
&\rho(x,0)=\rho^{0}(x),  &\mbox{in}&&I,\\
&\rho(0,t)=\rho(1,t)=0, &&&\forall t\in [0,\infty),
\end{aligned}
\right.
\end{equation}
and
\begin{equation}\label{hjbound}
\left\{
\begin{aligned}
&\kappa_{t}\kappa_{x}=\rho_{t}\rho_{x},\quad&\mbox{in}&&\quad 
I\times(0,T)&,\\
&\kappa(x,0)=\kappa^{0}(x),\quad &\mbox{in}&&\quad I,&\\
&\kappa(0,t)=\kappa(0,0)\quad\mbox{and}&\kappa&(1,t)=\kappa(1,0),
&\quad \forall t\in [0,T].
\end{aligned}
\right.
\end{equation}
Denote $I_T$ by: 
$$I_T=I\times(0,T).$$ 
There are two 
natural assumptions concerning $\rho^0$ and
$\kappa^0$, the first one is again the positivity of the dislocation
densities $\theta^+$ and $\theta^-$ at the initial time, which
yields to the following condition:
\begin{equation}\label{condrev}
\kappa^{0}_{x}\geq|\rho^{0}_{x}|,
\end{equation}
and the second one has to do with the balance of the physical model that 
starts
with the same number of positive and negative dislocations. In other
words, if $n^+$ and $n^-$ are the total number of positive and
negative dislocations respectively at $t=0$ then:
\begin{eqnarray*}
\rho^{0}(1)-\rho^{0}(0)&=&\int_{0}^{1}\rho^{0}_{x}(x)\,dx,\\
&=&\int_{0}^{1}(\theta^{+}(x,0)-\theta^{-}(x,0))\,dx,\\
&=& n^+-n^- = 0,
\end{eqnarray*}
this shows that  $\rho^{0}(1)=\rho^{0}(0)$ and this is what appears
in (\ref{heatbound}). Up to now, formal relations between the initial
conditions are only expressed. Whereas, required regularity,  together
with  the announcement of the main results will be
stated in the next subsection.
\subsection{Statement of the main results on a bounded interval}
From now on, the reader should not be confused with the term $\rho$ that
will always be the unique solution of the classical heat
equation (\ref{heatbound}). 
The two main theorems that we are going to prove are:
\begin{theorem}{\bf (Existence and uniqueness of a viscosity 
solution)}\label{boreps}\\
Let $T>0$ and $\eps>0$ be two constants. Take $\kappa^{0}\in
Lip(I)$ and $\rho^{0}\in C^{\infty}_{0}(I)$ satisfying:
$$\kappa^{0}_{x}\geq G(\rho^{0}_{x})\quad\mbox{a.e. in}\quad I,  $$
where
$$G(x)=\sqrt{x^{2}+\eps^{2}},$$
then there exists a viscosity solution $\kappa\in
Lip(\bar{I}_{T})$ of (\ref{hjbound}), unique among those satisfying:
\begin{equation}\label{bossis} \kappa_{x}\geq G(\rho_{x})
\quad\mbox{a.e. in }\; \bar{I}_{T}.
\end{equation}
\end{theorem}
\begin{theorem}{\bf (Existence of a viscosity solution)}\label{borfinal}\\
Let $T>0$ and $\kappa^{0}\in Lip(I)$. Under the
condition (\ref{condrev}) satisfied a.e. in $I$, there exists a
viscosity solution $\kappa\in Lip(\bar{I}_{T})$ of (\ref{hjbound})
satisfying:
$$\kappa_{x}\geq |\rho_{x}|, \quad
\mbox{a.e. in }\;\; \bar{I}_{T}.$$
\end{theorem}
\subsection{Preliminary results}
Before proceeding with the proof of our theorems, we have to introduce some essential
tools that are the core of the "extension and restriction" method that we are going to use.\\\\
{\bf \underline{Extension of $\rho$ over $\R\times [0, T]$.}}\\
Consider the function $\hat{\rho}$ defined on $[0,2]\times[0,T]$
by
\begin{equation}\label{rhotilde}
\hat{\rho}(x,t)=\left\{ \begin{aligned} &\rho(x,t)
\hspace{0.27cm}\quad \mbox{if}
\quad (x,t)\in \bar{I}_T,\\
& -\rho(2-x,t)\quad\mbox{otherwise},
\end{aligned}
\right.
\end{equation}
this is just a $C^{1}$ antisymmetry of $\rho$ with respect to the line
$x=1$. The continuation of $\hat{\rho}$ to $\R\times [0, T]$ is made
by spatial periodicity of period $2$. A simple computation yields, for
$(x,t)\in (1,2)\times (0,T)$:
$$\hat{\rho}_{t}(x,t)=-\rho_{t}(2-x,t)\quad\mbox{and}\quad
\hat{\rho}_{xx}(x,t)=-\rho_{xx}(2-x,t),$$
and hence it is easy to verify that
$\hat{\rho}\,|_{[1, 2]\times[0, T]}$ solves (\ref{heatbound}) with
$I$ replaced with the interval $(1, 2)$ and $\rho^{0}$ replaced with its
symmetry with respect to the point $x=1$; the boundary conditions
are unchanged and the regularity of the initial condition is
conserved. To be more precise, we write down
some useful properties of $\hat{\rho}$.\\\\
{\bf \underline{Regularity properties of $\hat{\rho}$.}}\\
Let $r$ and $s$ are two positive integers such that $s\leq 2$. From the 
construction of $\hat{\rho}$ and the above discussion, we get the
following:
\begin{equation}\label{rhotildereg}
\begin{aligned}
&\mbox{i})\, \hat{\rho}_t \;\mbox{and} \;\hat{\rho}_x \;\mbox{are
in} \;C(\R\times[0,T]),\\
&\mbox{ii})\, \hat{\rho}=0 \;\mbox{on}\;
\;\mathbb{Z}\times [0,T],\\
&\mbox{iii})\, \hat{\rho}_{t}=\hat{\rho}_{xx}\; \mbox{on}\;
(\R\setminus\mathbb{Z})\times(0,T),\\
&\mbox{iv})\,
||\partial_{t}^{r}\partial_{x}^{s}\hat{\rho}(.,t)||_{L^{\infty}(\R)}\leq
C, \;\forall t\in[0,T],
\end{aligned}
\end{equation} 
Where $C$ is a certain constant and the limitation $s\leq 2$ comes from the
spatial antisymmetry. These conditions are valid thanks to the way of 
construction of the function $\hat{\rho}$ and to the maximum
principle of the solution of the heat equation on bounded domains
(see \cite{brez1, evans}).\\ 

 Let
\begin{equation}\label{g_hat}
\hat{g}(x,t)=-\hat{\rho}_{t}(x,t)\hat{\rho}_{x}(x,t).
\end{equation}
From the above discussion, it is worth noticing that this function
is a Lipschitz continuous function in the $x$-variable.

The following three lemmas will be used in
the proof of Theorem \ref{boreps}.
\begin{lemma}{\bf (Entropy sub-solution)}\label{2entlem}\\
The function $G(\hat{\rho}_{x})$ is an entropy sub-solution of
\begin{equation}\label{2enteqn}
\left\{
\begin{aligned}
&w_{t}+(\hat{g}f_{\eps}(w))_{x}=0,\quad&\mbox{in}&\quad Q_T,\\
&w(x,0)=w^{0}(x)\quad&\mbox{in}&\quad \R,
\end{aligned}
\right.
\end{equation}
where $f_{\eps}$ is given by (\ref{functionf}), and
$w^{0}(x)=G(\hat{\rho}_{x}(x,0))$.
\end{lemma}
\noindent {\bf Proof.} Similar to Lemma \ref{strsubsol}.
$\hfill{\Box}$
\begin{lemma}{\bf (Differentiability property)}\label{aelemma}\\
Let $u(x,t)$ be a differentiable function with respect to $(x,t)$ a.e. in 
$Q_T$.
Define the set $M$ by:
$$M=\left\{x\in \R;\;\; u \;\mbox{is differentiable a.e. in}\; \{x\}\times(0,T)\right\}, $$
then $M$ is dense in $\R$.
\end{lemma}
\noindent {\bf Proof.} Define  $\mathcal{L}^n$,
$n\in\mathbb{N}$ to be the Lebesgue $n$-dimensional measure. Let
$N\subset Q_T$ be the set defined by:
$$N=\left\{(x,t)\in Q_T;\;\;u\;\mbox{is not differentiable on 
}(x,t)\right\},$$
and let $\mathbb{I}_{N}$ be the characteristic function of the set
$N$. Since $\mathcal{L}^2(N)=0$, we can write,
$$\int_{Q_T}\mathbb{I}_{N}(x,t)dxdt=0.$$
Using Fubini's theorem we get
$$\int_{\R}g(x)dx=0, \quad\mbox{with}\quad
g(x)=\left(\int_{0}^{T}\mathbb{I}_{N}(x,t)dt\right)\geq 0,$$ then
$$g=0 \quad\mbox{a.e. in } \quad\R $$
and consequently
$$J=\{x;\,\,g(x)\neq 0\} \quad\mbox{verifies}\quad\mathcal{L}^{1}(J)=0.$$
In other words,
$$\forall x\in \R\setminus J,\,\,\,u(x,\cdot)\,\,\mbox{is
   differentiable with respect to $(x,t)$ a.e. in}\,\,
(0,T),$$
hence $\R\setminus J\subset M$ which implies our lemma. $\hfill{\Box}$\\

In the next lemma, we show a lower-bound estimate for the gradient of
$\hat{\kappa}$ analogue to (\ref{bossis}). This was previously done for
$\kappa_{x}$ in the case where $g$ is a twice continuously differentiable
function using mainly Theorems \ref{relation} and \ref{entcomp}. Here, the way of extending the
function $\rho$ over $\bar{Q}_T$ makes $\hat{g}$ loose some
of the regularity stated in Theorem \ref{relation}. However, the
following lemma shows that a similar result holds in the case
$\hat{g}\in W^{1,\infty}(\bar{Q}_T)$.
\begin{lemma}\label{relation_bis}
The function $\hat{\kappa}_{x}\in L^{\infty}(Q_T)$ is an entropy solution of
(\ref{2enteqn}) with initial data $w^{0}=\hat{\kappa}^{0}_{x}\in
L^{\infty}(\R)$.
\end{lemma}
\noindent {\bf Proof of Lemma \ref{relation_bis}.} Let $\tilde{g}$ be an extension of the
function $\hat{g}$ on $\R^{2}$ defined by:
\begin{equation}\label{g_tilde}
\tilde{g}(x,t)=\left\{
\begin{aligned}
&\hat{g}(x,t)\quad&\mbox{if}&\quad(x,t)\in \bar{Q}_T,\\
&\hat{g}(x,T)\quad&\mbox{if}&\quad t>T,\\
&\hat{g}(x,0)\quad&\mbox{if}&\quad  t<0.
\end{aligned}
\right.
\end{equation}
Consider a sequence of mollifiers $\xi^{n}$ in $\R^{2}$ and let
$\tilde{g}^{n}=\tilde{g}*\xi^{n}$. Remark that, from the standard
properties of the mollifier sequence, we have
$\tilde{g}^{n}\in C^{\infty}(\R^2)$ and:
\begin{equation}\label{conver_tilde0}
\tilde{g}^{n}\rightarrow \hat{g}\;\;\mbox{uniformly on compacts in
}\;\bar{Q}_T,
\end{equation}
and\begin{equation}\label{conver_tilde1}
\tilde{g}^{n}_{x}\rightarrow \hat{g}_{x}\;\;\mbox{in
}\;L^{p}_{loc}(Q_T),\quad 1\leq p<\infty,
\end{equation}
together with the following estimates:
\begin{equation}\label{estim_tilde}
||\partial^{r}_{t}\partial^{s}_{x}\tilde{g}^{n}||_{L^{\infty}(\bar{Q}_T)}\leq
||\partial^{r}_{t}\partial^{s}_{x}\hat{g}||_{L^{\infty}(\bar{Q}_T)}\;\;\mbox{for}\;\;r,
s \in \mathbb{N},\;r+s\leq 1.
\end{equation}
Now, take again the Hamilton-Jacobi equation (\ref{hjextension}) with
$\hat{g}$ replaced with $\tilde{g}^{n}$:
\begin{equation}\label{hjextensiongg}
\left\{
\begin{aligned}
&u_{t}+\tilde{g}^{n}f_{\eps}(u_{x})=0\quad&\mbox{in}
&\quad \R\times (0,T),\\
&u(x,0)=\hat{\kappa}^{0}(x)\quad&\mbox{in}&\quad\R,
\end{aligned}
\right.
\end{equation}
and notice that the above properties of the function $\tilde{g}^{n}$
enters us into the framework 
of Theorem \ref{relation}. Thus, we have a unique viscosity solution
$\tilde{\kappa}^{n}\in Lip(\bar{Q}_T)$ of (\ref{hjextensiongg}) with initial
condition $\hat{\kappa}^{0}$ whose spatial
derivative $\tilde{\kappa}^{n}_{x}\in L^{\infty}(Q_T)$ is an entropy
solution of the corresponding derived equation with initial data
$\hat{\kappa}^{0}_{x}$. From Remark \ref{uniform_bound} and
(\ref{estim_tilde}), we deduce that the sequence
$(\tilde{\kappa}^{n})_{n\geq 1}$ is locally uniformly bounded in
$W^{1,\infty}(\bar{Q}_T)$ and that:
\begin{equation}\label{ktnx_bound}
||\tilde{\kappa}^{n}_{x}||_{L^{\infty}(Q_T)}\leq
||\hat{\kappa}^{0}_{x}||_{L^{\infty}(\R)}+T||\hat{g}_{x}||_{L^{\infty}(Q_T)}||f_{\eps}||_{L^{\infty}(\R)}.
\end{equation}
Moreover, from (\ref{conver_tilde0}), we use again the Stability Theorem
of viscosity solutions \cite[Theorem 2.3]{barbook}, and we obtain:
\begin{equation}\label{conv_ktn}
\tilde{\kappa}^{n}\rightarrow \hat{\kappa}\;\mbox{ locally uniformly in
}\; \bar{Q}_T.
\end{equation}
Back to the entropy solution, we write down the entropy inequality (see
Definition \ref{def_ent}) satisfied by $\tilde{\kappa}^{n}_{x}$:
\begin{equation}\label{eqn_ktnx}
\int_{Q_T}\Big(\,\eta(\tilde{\kappa}^{n}_{x})\phi_{t}+\Phi(\tilde{\kappa}^{n}_{x})\tilde{g}^{n}\phi_{x}
+h(\tilde{\kappa}^{n}_{x})\tilde{g}^{n}_{x}\phi\,\Big)dxdt+\int_{\R}\eta(\hat{\kappa}^{0}_{x})\phi(x,0)
dx\geq 0,
\end{equation}
where $\eta$, $\Phi$, $h$ and $\phi$ are given by Definition
\ref{def_ent}. Taking (\ref{ktnx_bound}) into consideration, we use
a property of bounded sequences in $L^{\infty}(Q_T)$ (see
\cite[Proposition 3]{egh2}) that guarantees the existence of a
subsequence (call it again $\tilde{\kappa}^{n}_{x}$) so that, for
any function $\psi\in C(\R;\R)$,
\begin{equation}\label{conpsi}
\psi(\tilde{\kappa}^{n}_{x})\rightarrow U_{\psi}\;\mbox{ weak$-\star$}\;\;\;\mbox{in}
\;\;\;L^{\infty}(Q_T).
\end{equation}
Furthermore, there exists $\mu\in L^{\infty}(Q_{T}\times(0,1))$
such that:
\begin{equation}\label{conmu}
\int_{0}^{1}\psi(\mu(x,t,\alpha))d\alpha=U_{\psi}(x,t),\;\mbox{ for
  a.e. }\;(x,t) \in Q_T.
\end{equation}
Applying (\ref{conpsi}) with $\psi$ replaced with $\eta$, $\Phi$ and
$h$ respectively, and using (\ref{conmu}), we get:
\begin{equation}\label{coneph}
\left\{
\begin{aligned}
&\eta(\tilde{\kappa}^{n}_{x}(.))\rightarrow
\int_{0}^{1}\eta(\mu(.,\alpha))d\alpha&\mbox{ weak$-\star$}&\;\;\;&\mbox{in}&\;\;\;L^{\infty}(Q_T),\\
&\Phi(\tilde{\kappa}^{n}_{x}(.))\rightarrow
\int_{0}^{1}\Phi(\mu(.,\alpha))d\alpha&\mbox{ weak$-\star$}&\;\;\;&\mbox{in}&\;\;\;L^{\infty}(Q_T),\\
&h(\tilde{\kappa}^{n}_{x}(.))\rightarrow
\int_{0}^{1}h(\mu(.,\alpha))d\alpha&\mbox{ weak$-\star$}&\;\;\;&\mbox{in}&\;\;\;L^{\infty}(Q_T).
\end{aligned}
\right.
\end{equation}
This, together with (\ref{conver_tilde0}), (\ref{conver_tilde1})
permits to pass to the limit in (\ref{eqn_ktnx}) in the
distributional sense, hence we get:
\begin{equation}\label{aft_con}
\begin{aligned}
\int_{Q_T}\int_{0}^{1}\Big(\eta(\mu(.,\alpha))\phi_{t}+\Phi(\mu(.,\alpha))\hat{g}\phi_{x}
+h(\mu(.,\alpha))\hat{g}_{x}\phi\Big)\,dxdtd\alpha+\\
\int_{\R}\eta(\hat{\kappa}^{0}_{x})\phi(x,0)
dx\geq 0.
\end{aligned}
\end{equation}
In \cite[Theorem 3]{egh2}, the function $\mu$ satisfying
(\ref{aft_con}) is called an entropy process solution. It has been proved to
be unique and independent of $\alpha$. Although this result in
\cite{egh2} was for a divergence-free function $\hat{g}\in
C^{1}(\bar{Q}_T)$, we remark that it can be adapted to the case of any
function $\hat{g}\in W^{1,\infty}(\bar{Q}_T)$ (see for instance
Remark \ref{final_rem} and the proof of \cite[Theorem 3]{egh2}).
Using this, we infer the existence of a function $z\in
L^{\infty}(Q_T)$ such that:
\begin{equation}\label{mu_ind_a}
z(x,t)=\mu(x,t,\alpha), \;\mbox{ for a.e. }\; (x,t,\alpha)\in
Q_{T}\times(0,1),
\end{equation}
hence, $z$ is an entropy solution of (\ref{2enteqn}). We now make use of
(\ref{mu_ind_a}) and we apply equality (\ref{conmu})
for $\psi(x)=x$ to obtain,
\begin{equation}\label{mu=ka*}
z=\,\mbox{weak$-\star$}\lim_{n\rightarrow
  \infty}\tilde{\kappa}^{n}_{x}\;\;\mbox{ in }\;\;L^{\infty}(Q_T).
\end{equation}
From (\ref{mu=ka*}) and (\ref{conv_ktn}) we deduce that,
$$z(x,t)=\hat{\kappa}_{x}(x,t)\;\;\mbox{ a.e. in }\;\;Q_{T},$$
which completes the proof of Lemma \ref{relation_bis}. $\hfill{\Box}$
\subsection{Proofs of Theorems \ref{boreps}, \ref{borfinal}}
\noindent {\bf Proof of Theorem \ref{boreps}.} We extend
the function $\kappa^{0}$ to $\hat{\kappa}^{0}\in Lip(\R)$
in the following way:
\begin{equation}\label{newkapa}
\hat{\kappa}^{0}(x)=\left\{
\begin{aligned}
&\kappa^{0}(x)\quad &\mbox{if}&\quad x\in[0,1],\\
&(||\rho^{0}_{x}||_{L^{\infty}(I)}+\eps)(x-1)+\kappa^{0}(1)\quad&\mbox{if}&\quad
x\geq1,\\
&(||\rho^{0}_{x}||_{L^{\infty}(I)}+\eps)x+\kappa^{0}(0)\quad&\mbox{if}&\quad
x\leq0.
\end{aligned}
\right.
\end{equation}
Consider the initial value problem defined by:
\begin{equation}\label{hjextension}
\left\{
\begin{aligned}
&u_{t}+\hat{g}f_{\eps}(u_{x})=0\quad&\mbox{in}
&\quad \R\times (0,T),\\
&u(x,0)=\hat{\kappa}^{0}(x)\quad&\mbox{in}&\quad\R.
\end{aligned}
\right.
\end{equation}
This is a Hamilton-Jacobi equation with a Hamiltonian $F\in
C(\bar{Q}_{T}\times \R)$ defined by:
$$F(x,t,u)=\hat{g}(x,t)f_{\eps}(u).$$
From the regularity properties of $\hat{\rho}$, we can directly see
that $\bold{(V0)}$-$\bold{(V1)}$-$\bold{(V2)}$ are satisfied; this
is quite similar to what was done in Proposition \ref{approblem}.
Since $\hat{\kappa}^{0}$ is a Lipschitz continuous function, we deduce from Theorems
\ref{comphj}, \ref{exishj} and Proposition \ref{lipreghj} the existence and uniqueness of a
viscosity solution $\hat{\kappa}\in Lip(\bar{Q}_T)$ of
(\ref{hjextension}). Moreover, in order to recover the boundary
conditions given by (\ref{hjbound}) on $\partial I\times[0,T]$, we proceed as follows. Let $M$ be the set
defined by Lemma \ref{aelemma} and let $x\in M$. For every $t\in[0,T]$, we
write:
\begin{equation}
|\hat{\kappa}(x,t)-\hat{\kappa}(x,0)|\leq\int_{0}^{t}|\hat{\kappa}_{s}(x,s)|ds
\leq\int_{0}^{t}|F(x,s,\hat{\kappa}_{x}(x,s))|ds\leq
\int_{0}^{t}\left(|F(0,s,\hat{\kappa}_{x}(x,s))|+C|x|\right)ds.
\nonumber
\end{equation}
In these inequalities we have used the fact that $\hat{\kappa}$ is a
Lipschitz continuous viscosity solution of (\ref{hjextension}) and hence it
verifies the equation in $Q_T$ at the points where it is
differentiable (see for instance \cite{barbook}). Also, we have used
the condition $\bold{(F1)}$ with $p=q$ and $C_{R}=C$, a constant
independant of $R$. Now from (\ref{rhotildereg})-(ii), we deduce that:
$$|F(0,s,\hat{\kappa}_{x}(x,s))|=|\hat{\rho}_{x}(0,s)\hat{\rho}_{t}(0,s)f_{\eps}(\hat{\kappa}_{x}(x,s))|=0,
\quad \mbox{for a.e.}\quad s\in(0,t),$$ and hence we get
\begin{equation}\label{bo7a}
|\hat{\kappa}(x,t)-\hat{\kappa}(x,0)|\leq C|x|t.
\end{equation}
Since $M$ is a dense subset of $\R$, we pass to the limit in
(\ref{bo7a}) as $x\rightarrow0$ and the equality
$$\hat{\kappa}(0,t)=\hat{\kappa}(0,0)=\kappa^{0}(0) \quad\forall t\in[0,T]$$
holds. Similarly, we can verify that $\hat{\kappa}(1,t)=
\hat{\kappa}(1,0)=\kappa^{0}(1)$ for all $t\in[0,T]$.\\ 

\noindent {\bf\underline{Existence.} }The extension $\hat{\kappa}^{0}$ of $\kappa^{0}$ outside the
interval $I$ is a linear extension of slope~$||\rho^{0}_{x}||_{L^{\infty}(I)}+~\eps$, therefore we have,
\begin{equation}\label{initcond2}
\hat{\kappa}^{0}_{x}(\cdot)\geq
\sqrt{(\hat{\rho}^{0}_{x}(\cdot))^{2}+\eps^{2}}=G(\hat{\rho}^{0}_{x}(\cdot)),
\quad \mbox{a.e. in }\R.
\end{equation}
From Lemma \ref{relation_bis}, we know that $\hat{\kappa}_{x}$ is an
entropy solution of equation
(\ref{2enteqn}) and from Lemma \ref{2entlem}, we know that
$G(\hat{\rho}_{x})$ is an entropy
sub-solution of (\ref{2enteqn}). Since (\ref{initcond2}) holds, we use
the Comparison Theorem \ref{entcomp} to get,
\begin{equation}\label{akbarmin}
\hat{\kappa}_{x}(x,t)\geq
\sqrt{\hat{\rho}^{2}_{x}(x,t)+\eps^{2}}\geq\eps>0, \quad \mbox{ for
a.e. } (x,t)\in \bar{Q}_T.
\end{equation}
\noindent Take $\kappa $ to be the restriction of $\hat{\kappa}$ on
$\bar{I}_{T}$ where $\hat{\kappa}^{0}$ and $\hat{\rho}$ have their
automatic replacements $\kappa^{0}$ and $\rho$ respectively on this
subdomain. It is clear that $\kappa\in Lip(\bar{I}_{T})$ is a viscosity
solution of:
\begin{equation}\label{wneheyto}
\left\{
\begin{aligned}
&\kappa_{t}+gf_{\eps}(\kappa_{x})=0\quad&\mbox{in}&\quad I_T,\\
&\kappa(x,0)=\kappa^{0}(x)\quad&\mbox{in}&\quad I,\\
&\kappa(0,t)=\kappa^{0}(0)\quad\mbox{and}\quad
\kappa(1,t)=\kappa^{0}(1)\;\;&\forall&\;0\leq t \leq T,
\end{aligned}
\right.
\end{equation}
where $g(x,t)=-\rho_{t}(x,t)\rho_{x}(x,t)$ and $\kappa_{x}(x,t)\geq
G(\rho_{x}(x,t))$ for a.e. $(x,t)\in \bar{I}_T$. We also notice that
$\kappa$ is a viscosity solution of (\ref{hjbound}), for it suffices
to follow the same steps of the passage from the viscosity solution
of (\ref{realappeqn}) to the viscosity solution of
(\ref{mainhj}) (see the proof of Theorem \ref{main1} for details).\\

\noindent{\bf \underline{Uniqueness.}} Since the function
$$\bar{H}(x,t,u)=g(x,t)f_{\eps}(u)\in
C(\bar{I}_{T}\times\R)$$
satisfies for a fixed $t$: $$|\bar{H}(x,t,u)-\bar{H}(y,t,u)|\leq C(|x-y|(1+|u|)),$$
for every $x,\;y\in (0,1)$ and $u\in\R$, we use \cite[Theorem 
2.8]{barbook} to show that $\kappa$ is the unique viscosity solution of
(\ref{wneheyto}). We claim that $\kappa$ is the unique viscosity solution of
(\ref{hjbound}). Indeed, we can also follow the same
mechanism as in the proof of Theorem \ref{main1}.$\hfill{\Box}$\\

We now move towards the proof of Theorem \ref{borfinal} that has the same 
flavor of what was done in Section \ref{sec4}. We just need to care about the change in the
structure of our problem and the boundary conditions. Our first step
will be the following lemma.
\begin{lemma}\label{tekrar}
Let $c_1$ and $c_2$ be two positive constants defined respectively
by:
$$c_{1}=||(\rho^{0}_{xx})^{2}||_{L^{\infty}(I)}+||\rho^{0}_{x}||_{L^{\infty}(I)}||\rho^{0}_{xxx}||_{L^{\infty}(I)},$$
and
$$c_{2}=(||\kappa^{0}_{x}||_{L^{\infty}(I)}+1)^{2}.$$
Then the function $\bar{S}$ defined on $Q_T$ by:
$$\bar{S}(x,t)=\sqrt{2c_{1}t+c_{2}}$$
is an entropy super-solution of (\ref{2enteqn})
with $$w^{0}(x)=\bar{S}(x,0)=||\kappa^{0}_{x}||_{L^{\infty}(I)}+1.$$
\end{lemma}
\noindent {\bf Proof.} See Subsection 5.1-III.$\hfill{\Box}$\\

\noindent {\bf Proof of Theorem \ref{borfinal}.} Let
$\eps>0$ be a fixed constant. Define $\hat{\kappa}^{0,\eps}\in
Lip(\R)$ by:
\begin{equation}\label{newkapaeps}
\hat{\kappa}^{0,\eps}(x)=\left\{
\begin{aligned}
&\kappa^{0}(x)+\eps x\quad&\mbox{if}&\quad x\in[0,1],\\
&(||\kappa^{0}_{x}||_{L^{\infty}(I)}+\eps)(x-1)+(\kappa^{0}(1)+\eps)\quad&\mbox{if}&\quad 
x\geq1,\\
&(||\kappa^{0}_{x}||_{L^{\infty}(I)}+\eps)x+\kappa^{0}(0)\quad&\mbox{if}&\quad 
x\leq0.
\end{aligned}
\right.
\end{equation}
Since $\kappa^{0}_{x}\geq |\rho^{0}_{x}|\; \mbox{a.e. }\mbox{in }I$,
it is clear that for a.e. $x\in \R$ we have
$$\hat{\kappa}^{0,\eps}_{x}\geq G(\hat{\rho}^{0}_{x}),$$
and hence, from the discussion of the proof of Theorem \ref{boreps},
there exists a unique viscosity solution $\hat{\kappa}^{\eps}\in
Lip(\bar{Q}_T)$ of 
\begin{equation}\label{wnehayto} \left\{
\begin{aligned}
&\hat{\kappa}^{\eps}_{t}\hat{\kappa}^{\eps}_{x}=\hat{\rho}_{t}\hat{\rho}_{x}\quad&\mbox{in}&\quad Q_T,\\
&\hat{\kappa}^{\eps}(x,0)=\hat{\kappa}^{0,\eps}(x)\in Lip(\mathbb{R})\quad&\mbox{in}&\quad\mathbb{R},
\end{aligned}
\right.
\end{equation}
unique among those satisfying:
\begin{equation}\label{7amoudiel7akim}
\hat{\kappa}^{\eps}_{x}\geq G(\hat{\rho}_{x})\;\,\;\;\mbox{a.e.
}\mbox{in } \;\;\bar{Q}_{T}.
\end{equation}
Assume without loss of generality that $\eps<1$. The $\eps$-uniform bound 
for
$\hat{\kappa}^{\eps}_{t}$ is trivial, it suffices to use directly
the equation satisfied by $\hat{\kappa}^{\eps}$  together with
(\ref{7amoudiel7akim}). And the $\eps$-uniform bound for
$\hat{\kappa}^{\eps}_{x}$ follows from Lemma \ref{tekrar} and
Theorem \ref{entcomp} since
$$\hat{\kappa}^{\eps}_{x}(x,0)\leq
||\kappa^{0}_{x}||_{L^{\infty}(I)}+\eps\leq||\kappa^{0}_{x}||_{L^{\infty}(I)}+1=
\sqrt{c_2}=\bar{S}(x,0).$$ Following exactly the same technic of Section
\ref{sec4}, namely the proof of Theorem \ref{maintheorem}, we get that the
sequence $\hat{\kappa}^{\eps}$ converges locally uniformly to
$\hat{\kappa}$ in $\bar{Q}_{T}$ with $\hat{\kappa}\in
Lip(\bar{Q}_{T})$ satisfies,
\begin{equation}\label{tarre1}
\hat{\kappa}_{x}\geq |\hat{\rho}_{x}|\;\;\;\mbox{a.e.}\mbox{ in
}\;\; \bar{Q}_{T}
\end{equation}
and
\begin{equation}\label{tarre2}
\hat{\kappa}(x,0)=\hat{\kappa}_{0}(x)\;\; \mbox{ in }\;\;\R,
\end{equation} where $\hat{\kappa}_{0}$ is the uniform limit
of the sequence $\hat{\kappa}^{0,\eps}$ in $\R$. Theorem
\ref{boreps} guarantees that
\begin{equation}\label{a5erboreps}
\hat{\kappa}^{\eps}(0,t)=\hat{\kappa}^{0,\eps}(0)=\kappa^{0}(0),
\end{equation}
and
\begin{equation}\label{a5erboreps1}
\hat{\kappa}^{\eps}(1,t)=\hat{\kappa}^{0,\eps}(1)=\kappa^{0}(1)+\eps,
\end{equation}
for all $t\in [0,T]$. From (\ref{a5erboreps}), (\ref{a5erboreps1})
and the pointwise convergence, up to a subsequence, of $\hat{\kappa}^{\eps}$ 
to
$\hat{\kappa}$, we deduce that
\begin{equation}\label{a5erboreps2}
\hat{\kappa}(0,t)=\lim_{\eps\rightarrow0}\hat{\kappa}^{\eps}(0,t)=\kappa^{0}(0),\quad\forall
t\in [0,T],
\end{equation}
and
\begin{equation}\label{a5erboreps3}
\hat{\kappa}(1,t)=\lim_{\eps\rightarrow0}\hat{\kappa}^{\eps}(1,t)=
\lim_{\eps\rightarrow0}(\kappa^{0}(1)+\eps)=\kappa^{0}(1)
\quad\forall t\in [0,T].
\end{equation}
Take $\kappa$ to be the restriction of $\hat{\kappa}$ over
$\bar{I}_T$; $\hat{\rho}$ and $\hat{\kappa}_{0}$ have their
automatic replacements $\rho$ and $\kappa^{0}$ respectively on this
restricted domain. From (\ref{tarre1}), (\ref{tarre2}),
(\ref{a5erboreps2}) and (\ref{a5erboreps3}), we deduce that
$\kappa$ is the required solution.$\hfill{\Box}$
\end{section}
\begin{section}{Appendix: Proof of Theorem \ref{entcomp}}\label{sec6}
We will work on the entropy inequality (\ref{withsgn}) satisfied by
$u$ and its analogue satisfied by $v$, using the dedoubling variable technique 
of Kruzhkov (see
\cite{kru1}) and following the same steps of \cite[Theorem 3]{egh2},
taking into consideration the new modifications arising from the fact
that we are dealing with sub-/super-entropy solutions and the fact that
$g\in W^{1,\infty}(\bar{Q}_T)$ is not a gradient-free function.\\
\noindent The proof can be divided into three steps. Denote $B_r$ by 
$B_{r}=\{x\in \R; \;\;|x|\leq r\}$ for any $r>0$,
$F^{\pm}(u,v)=sgn^{\pm}(u-v)(f(u)-f(v))$,
\begin{equation}\label{max_function}
y^{\infty}=||y||_{L^{\infty}(Q_T)}\quad\mbox{ for every }\quad y\in 
L^{\infty}(Q_T)
\end{equation}
and
\begin{equation}\label{max_f}
M_{f}=\max_{|x|\leq
  \max(u^{\infty},v^{\infty})}|f^{'}(x)|.
\end{equation}
In step 1, we prove that the initial conditions $u^{0}$, $v^{0}$ satisfy
for any $a>0$:
\begin{equation}\label{limfor_u0}
\lim_{\tau\rightarrow 
0}\frac{1}{\tau}\int_{0}^{\tau}\int_{B_{a}}(u(x,t)-u^{0}(x))^{+}dxdt=0,
\end{equation}
\begin{equation}\label{limfor_v0}
\lim_{\tau\rightarrow 
0}\frac{1}{\tau}\int_{0}^{\tau}\int_{B_{a}}(v(x,t)-v^{0}(x))^{-}dxdt=0,
\end{equation}
respectively.\\\\
\noindent In step 2, The following relation between $u$ and $v$ is
shown:
\begin{equation}\label{relation_u_v}
\int_{Q_T}\left[(u(x,t)-v(x,t))^{+}\psi_{t}+F^{+}(u(x,t),v(x,t))g(x,t)\psi_{x}\right]dxdt\geq
0,
\end{equation}
for every $\psi\in C^{1}_{0}(\R\times(0,T);\R_{+})$.\\
\noindent After that, we define $A(t)$ for $0<t<\min(T,\frac{a}{\omega})$ 
and
$\omega=g^{\infty}M_{f}$, by:
\begin{equation}\label{A_t}
A(t)=\int_{B_{a-{\omega t}}}\left(u(x,t)-v(x,t)\right)^{+}dx.
\end{equation}
In step 3, we show that $A$ is non-increasing a.e. in
$(0,\min(T,\frac{a}{\omega}))$ and we deduce that
$$u(x,t)\leq v(x,t)\quad \mbox{a.e. in} \quad Q_T.$$\\
{\bf Step 1: Proof of (\ref{limfor_u0}), (\ref{limfor_v0}).}\\\\
Let $\xi^{n}$ be a sequence of mollifiers in $\R$ with $\xi^{1}=\xi$.
Recall that the function $\xi\in C^{\infty}_{0}(\R)$ satisfies the
following properties:
\begin{equation}\label{propmol}
\begin{aligned}
&\hspace{-4cm}\mbox{supp}(\xi)=\{x\in \R,\;\xi(x)\neq 0\}\subset B_{1};\\
&\hspace{-4cm}\xi \geq 0,\;\quad\xi(-x)=\xi(x);\\
&\hspace{-4cm}\int_{B_{1}}\xi(x)dx=1;\\
&\hspace{-4cm}\xi^{n}(x)=n\xi(nx).
\end{aligned}
\end{equation}
Let $\tau\in \R$
such that $0<\tau<T$ and define the function $\gamma$ by:
\begin{equation}\label{gamma}
\gamma(t)=
\left\{
\begin{aligned}
&\frac{\tau-t}{\tau}\quad&\mbox{if}&\quad 0\leq t\leq \tau,\\
&0\quad&\mbox{if}&\quad t>\tau.
\end{aligned}
\right.
\end{equation}
Take $a>0$ and a test function $\psi\in C^{\infty}_{0}(\R;\R_{+})$ such
that,
$$\psi(x)=1\quad\mbox{for}\quad x\in B_{a}.$$
Let $y\in \R$ be a Lebesgue point of $u^{0}$ and we make
use of inequality (\ref{withsgn}) with $k=u^{0}(y)$ and the test
function $\phi(x,t)=\psi(x)\gamma(t)\xi^{n}(x-y)$ (this is possible since
$\phi$ is a permissible test function). Integrating the resulting
inequality with respect to $y$ over $\R$ yields:
\begin{equation}\label{sumT1234}
{\T}_{1}(n,\tau)+{\T}_{2}(n,\tau)+{\T}_{3}(n,\tau)+{\T}_{4}(n)\geq 0,
\end{equation}
with
\begin{equation}\label{tau1}
{\T}_{1}(n,\tau)=-\frac{1}{\tau}\int_{0}^{\tau}\int_{\R^2}(u(x,t)-u^{0}(y))^{+}\psi(x)\xi^{n}(x-y)\,dxdydt,
\end{equation}
\begin{equation}\label{tau2}
{\T}_{2}(n,\tau)=\int_{0}^{\tau}\int_{\R^2}F^{+}(u(x,t),u^{0}(y))g(x,t)\gamma(t)(\psi(x)\xi^{n}(x-y))_{x}dxdydt,
\end{equation}
\begin{equation}
\hspace{-3.8cm}{\T}_{3}(n,\tau)=-\int_{0}^{\tau}\int_{\R^2}sgn^{+}(u(x,t)-u^{0}(y))f(u^{0}(y))
\nonumber
\end{equation}
\begin{equation}\label{tau3}
\hspace{6cm}g_{x}(x,t)\gamma(t)\psi(x)\xi^{n}(x-y)dxdydt
\end{equation}
and
\begin{equation}\label{tau4}
{\T}_{4}(n)=\int_{\R^2}(u^{0}(x)-u^{0}(y))^{+}\psi(x)\xi^{n}(x-y)dxdy.
\end{equation}
Using the change of variables: $x=x^{'}$, $y=x^{'}-\frac{y^{'}}{n}$ in
(\ref{tau1}), and denoting again by $(x,y)$ the new variables
$(x^{'},y^{'})$ yields:
\begin{equation}\label{newtau1}
{\T}_{1}(n,\tau)=-\frac{1}{\tau}\int_{0}^{\tau}\int_{B_{1}}\int_{\R}\left(u(x,t)-u^{0}
\left(x-\frac{y}{n}\right)\right)^{+}
\psi(x)\xi(y)\,dxdydt,
\end{equation}
Using that,
\begin{equation}\label{strinq+}
(u-v)^{+}-(u-w)^{+}\leq (w-v)^{+} \quad \forall u, v, w\in \R,
\end{equation}
we infer that:
\begin{equation}
\hspace{-3.5cm}\T_{1}(n,\tau)+\overbrace{\frac{1}{\tau}\int_{0}^{\tau}\int_{\R}
(u(x,t)-u^{0}(x))^{+}\psi(x)dxdt}^{\T^{*}(\tau)}\leq
\nonumber
\end{equation}
\begin{equation}\label{tau1+}
\hspace{3.7cm} \psi^{\infty}\int_{K_{\psi}}\int_{B_1}\left|u^{0}
\left(x-\frac{y}{n}\right)-u^{0}(x)\right|\xi(y)dydx,
\end{equation}
where $K_{\psi}$ is the support of $\psi$. Same upper-bound,
independent of $\tau$, could be obtained for $\T_{4}(n)$.
Furthermore, since $u^{0}\in L^{\infty}(\R)$, thus integrable over
$K_{\psi}$, we use the Lebesgue differentiation Theorem to show that
the right side of (\ref{tau1+}) tends to $0$ when $n$ becomes large.
Now, let $\eps>0$, $\exists n_{0}$ such that
\begin{equation}\label{tawet14}
\T_{1}(n_{0},\tau)+\T^{*}(\tau)<\frac{\eps}{4}\quad\mbox{and}\quad
\T_{4}(n_{0})<\frac{\eps}{4},\;\;\;\;\forall \tau>0.
\end{equation}
We  also remark that the integrands of the right hand sides of
(\ref{tau2}) and (\ref{tau3}) are bounded and hence, for this particular
$n_{0}$ we can choose some $\tau_{0}$ such that
$\forall\;0<\tau<\tau_{0}$, we have:
\begin{equation}\label{tawet23}
\T_{2}(n_{0},\tau)<\frac{\eps}{4}\quad\mbox{and}\quad\T_{3}(n_{0},\tau)<\frac{\eps}{4}.
\end{equation}
From (\ref{tawet14}), (\ref{tawet23}) and (\ref{sumT1234}), we infer
that,
$$0<\T^{*}(\tau)< \eps,\quad\forall 0<\tau<\tau_{0}.$$
Since $\psi(x)=1$ over $B_{a}$, (\ref{limfor_u0}) is proven. Arguing in
the same way, we can prove (\ref{limfor_v0}). The slight difference is
using a similar inequality of (\ref{strinq+}) with $(\cdot)^{+}$
replaced with  $(\cdot)^{-}$.\\\\
{\bf Step 2: Proof of (\ref{relation_u_v}).}\\\\
It suffices to prove (\ref{relation_u_v}) for any function $\psi\in
C^{\infty}_{0}(Q_T;\R_+)$. We may also assume, without loss of
generality, that there is some $c>0$ such that $\psi(x,t)=0$ for
$t\in (0,c)\cup(T-c,T)$. For $n>\frac{1}{c}$, let $\xi^{n}$ be the
usual mollifier sequence in $\R$ and consider the function
$\phi(x,t,y,s)$ defined for $(x,t)\in Q_T$ and $(y,s)\in Q_T$ by,
$$\phi(x,t,y,s)=\psi\left(\frac{x+y}{2},\frac{t+s}{2}\right)\xi^{n}(x-y)\xi^{n}(t-s).$$
The function $\phi$ hence satisfies
$$\phi(.,.,y,s)\in C^{\infty}_{0}(Q_{T};\R_{+})\quad\mbox{and}
\quad\phi(x,t,.,.)\in C^{\infty}_{0}(Q_{T};\R_{+}).$$ Fix some
$(y,s)\in Q_T$ for which the function $v$ is well defined (this is
valid almost everywhere). Since $u$ is an entropy sub-solution of
(\ref{CL}), we consider the relation (\ref{withsgn}) satisfied by
$u$ with $k=v(y,s)$ and the test function $\phi(.,.,y,s)$. Upon
integrating this inequality with respect to $(y,s)$ over $Q_T$, we
get:
\begin{equation}
\int_{Q^{2}_{T}
}\left\{(u(x,t)-v(y,s))^{+}\phi_{t}(x,t,y,s)+F^{+}(u(x,t),v(y,s))
g(x,t)\phi_{x}(x,t,y,s)\right. \nonumber
\end{equation}
\begin{equation}\label{sub+}
\left.-sgn^{+}(u(x,t)-v(y,s))f(v(y,s))g_{x}(x,t)\phi(x,t,y,s)\right\}dxdtdyds
\geq 0.
\end{equation}
Similar inequality could be obtained since $v$ is an entropy
super-solution of (\ref{CL}). We just swap $+$, $u$ and $(x,t)$ with
$-$, $v$ and $(y,s)$ respectively, hence:
\begin{equation}
\int_{Q^{2}_T}\left\{(v(y,s)-u(x,t))^{-}\phi_{s}(x,t,y,s)+F^{-}(v(y,s),u(x,t))
g(y,s)\phi_{y}(x,t,y,s)\right. \nonumber
\end{equation}
\begin{equation}\label{sup-}
\left.\hspace{-0.05cm}-sgn^{-}(v(y,s)-u(x,t))f(u(x,t))g_{x}(y,s)\phi(x,t,y,s)\right\}dxdtdyds
\geq 0.
\end{equation}
Summing (\ref{sub+}) and (\ref{sup-}) and using the elementary
identities:
$$x^{-}=(-x)^{+}\quad\mbox{and}\quad sgn^{-}(x)=-sgn^{+}(-x), \quad\forall 
x\in \R,$$
we get, for $u=u(x,t)$ and $v=v(y,s)$,
\begin{equation}\label{after_sum}
\mathcal{Z}_{1}+\mathcal{Z}_{2}+\mathcal{Z}_{3}\geq 0,
\end{equation}
with:
\begin{equation}
\mathcal{Z}_{1}=\int_{Q^{2}_T}(u-v)^{+}(\phi_{t}+\phi_{s})(x,y,t,s)dxdtdyds,
\end{equation}
\begin{equation}
\mathcal{Z}_{2}=\int_{Q^{2}_T}F^{+}(u,v)[g(x,t)\phi_{x}(x,y,t,s)+g(y,s)\phi_{y}(x,y,t,s)]dxdtdyds,
\end{equation}
\begin{equation}
\mathcal{Z}_{3}=\int_{Q^{2}_T}sgn^{+}(u-v)[f(u)g_{x}(y,s)-f(v)g_{x}(x,t)]\phi(x,y,t,s)dxdtdyds.
\end{equation}
We now compute the first partial derivatives of the function $\phi$. For
$(x,t,y,s)\in Q_{T}\times Q_{T}$, we have:
\begin{equation}
\hspace{-3cm}\phi_{t}(x,t,y,s)=\xi^{n}(x-y)\left(\frac{1}{2}\psi_{t}\left(\frac{x+y}{2},\frac{t+s}{2}\right)\xi^{n}(t-s)\right.
\nonumber
\end{equation}
\begin{equation}
\hspace{5.7cm}\left.+\psi\left(\frac{x+y}{2},\frac{t+s}{2}\right)\xi^{n^{'}}(t-s) 
\right),
\end{equation}
\begin{equation}
\hspace{-3cm}\phi_{s}(x,t,y,s)=\xi^{n}(x-y)\left(\frac{1}{2}\psi_{t}\left(\frac{x+y}{2},
\frac{t+s}{2}\right)\xi^{n}(t-s)\right.
\nonumber
\end{equation}
\begin{equation}
\hspace{5.7cm}\left.-\psi\left(\frac{x+y}{2},\frac{t+s}{2}\right)\xi^{n^{'}}(t-s) 
\right),
\end{equation}
\begin{equation}
\hspace{-3cm}\phi_{x}(x,t,y,s)=\xi^{n}(t-s)\left(\frac{1}{2}\psi_{x}\left(\frac{x+y}{2},
\frac{t+s}{2}\right)\xi^{n}(x-y)\right.
\nonumber
\end{equation}
\begin{equation}
\hspace{5.7cm}\left.+\psi\left(\frac{x+y}{2},\frac{t+s}{2}\right)\xi^{n^{'}}(x-y) 
\right),
\end{equation}
\begin{equation}
\hspace{-3cm}\phi_{y}(x,t,y,s)=\xi^{n}(t-s)\left(\frac{1}{2}\psi_{x}\left(\frac{x+y}{2},
\frac{t+s}{2}\right)\xi^{n}(x-y)\right.
\nonumber
\end{equation}
\begin{equation}
\hspace{5.7cm}\left.-\psi\left(\frac{x+y}{2},\frac{t+s}{2}\right)\xi^{n^{'}}(x-y) 
\right).
\end{equation}
Using these relations in (\ref{after_sum}) and performing the following 
change of variables,
$$x^{'}=(x+y)/2,\;y^{'}=n(x-y),\;t^{'}=(t+s)/2,\;s^{'}=n(t-s);$$
denote the new variables $x^{'}$, $t^{'}$, $y^{'}$, $s^{'}$ by $x $, $t$,
$y$, $s$ and $\mathcal{Q}_{4}=Q_{T}\times B_{1}^{2}$. Also, for the
simplicity of expressions, denote
$$x^{+}=x+\frac{y}{2n},\;\;t^{+}=t+\frac{s}{2n},\;\;x^{-}=x-\frac{y}{2n},\;\;t^{-}=t-\frac{s}{2n}.$$
This altogether yields:
\begin{equation}\label{new_sum}
\mathcal{X}_{1}+\mathcal{X}_{2}+\mathcal{X}_{3}+\mathcal{X}_{4}\geq 0,
\end{equation}
with:
\begin{equation}\label{sum_new_1}
\mathcal{X}_{1}=\int_{\mathcal{Q}_{4}}(u(x^{+},t^{+})-v(x^{-},t^{-}))^{+}\psi_{t}(x,t)\xi(y)\xi(s)dxdtdyds,
\end{equation}
\begin{equation}\label{sum_new_2}
\begin{aligned}
\mathcal{X}_{2}=\frac{1}{2}\int_{\mathcal{Q}_{4}}F^{+}(u(x^{+},t^{+}),v(x^{-},t^{-}))(g(x^{+},t^{+})+g(x^{-},t^{-}))\times\\
\psi_{x}(x,t)\xi(y)\xi(s)dxdtdyds,
\end{aligned}
\end{equation}
\begin{equation}\label{sum_new_3}
\begin{aligned}
\mathcal{X}_{3}=\int_{\mathcal{Q}_{4}}F^{+}(u(x^{+},t^{+}),v(x^{-},t^{-}))(g(x^{+},t^{+})-g(x^{-},t^{-}))\times\\
\psi(x,t)n\xi^{'}(y)\xi(s)dxdtdyds,
\end{aligned}
\end{equation}
\begin{equation}\label{sum_new_4}
\begin{aligned}
\mathcal{X}_{4}=\int_{\mathcal{Q}_{4}}sgn^{+}(u(x^{+},t^{+})-v(x^{-},t^{-}))\left[f(u(x^{+},t^{+}))g_{x}(x^{-},t^{-})-\right.\\
\left.f(v(x^{-},t^{-}))g_{x}(x^{+},t^{+})\right]\psi(x,t)\xi(y)\xi(s)
dxdtdyds.
\end{aligned}
\end{equation}
At this point, it is worth mentioning that we will frequently use the 
following Lemma from \cite{kru2}.
\begin{lemma}\label{llip}
If $\Gamma\in Lip(\R)$ satisfies $|\Gamma(u)-\Gamma(v)|\leq C_{0}|u-v|$, then the function
$$H(u,v)=sgn^{+}(u-v)(\Gamma(u)-\Gamma(v))$$ satisfies
$|H(u,v)-H(u^{'},v^{'})|\leq C_{0}(|u-u^{'}|+|v-v^{'}|$ (see \cite[Lemma
3]{kru1}).
\end{lemma}
\noindent Consider now (\ref{sum_new_1}). Since 
$(u-v)^{+}=sgn^{+}(u-v)(u-v)$, we make
use of Lemma \ref{llip} to obtain:
\begin{equation}
\begin{aligned}
&\left|\mathcal{X}_{1}-\int_{Q_T}(u(x,t)-v(x,t))^{+}\psi_{t}(x,t)dxdt\right|\leq\\
&\left\{\int_{K_{\psi}}\int_{B_{1}^{2}}|u(x^{+},t^{+})-u(x,t)|(\psi_{t})^{\infty}\xi(y)\xi(s)dxdtdyds\right.\\
&\left.+\int_{K_{\psi}}\int_{B_{1}^{2}}|v(x^{-},t^{-})-v(x,t)|(\psi_{t})^{\infty}\xi(y)\xi(s)dxdtdyds\right\},
\end{aligned}
\nonumber
\end{equation}
where, by the Lebesgue Differentiation/Dominated Theorems, the right
hand side of this inequality tends to $0$ as $n \rightarrow \infty$,
and hence:
\begin{equation}\label{insh1}
\mathcal{X}_{1}\rightarrow 
\int_{Q_T}(u(x,t)-v(x,t))^{+}\psi_{t}(x,t)dxdt\quad\mbox{as}\;\;n\rightarrow\infty.
\end{equation}
Let us now turn to (\ref{sum_new_2}); using the fact that $g\in
W^{1,\infty}(Q_T)$ and hence Lipschitz continuous over the compact $K_{\psi}$, and
the fact that $F^{+}(u,v)$ is Lipschitz continuous in $u$ and $v$ (see Lemma
\ref{llip}), we get:
\begin{equation}
\begin{aligned}
&\left|\mathcal{X}_{2}-\int_{Q_T}F^{+}(u(x,t),v(x,t))g(x,t)\psi_{x}(x,t)dxdt
\right|\leq\\
&g^{\infty}M_{f}\psi_{x}^{\infty}
\left\{\int_{K_{\psi}}\int_{B_{1}^{2}}|u(x^{+},t^{+})-u(x,t)|\xi(y)\xi(s)dxdtdyds\right.\\
&\left.+\int_{K_{\psi}}\int_{B_{1}^{2}}|v(x^{-},t^{-})-v(x,t)|\xi(y)\xi(s)\right\}dxdtdyds\\
&+\frac{1}{n}C((g_{x})^{\infty},(g_{t})^{\infty},(\psi_{x})^{\infty},M_{f},u^{\infty},v^{\infty},T),
\end{aligned}
\end{equation}
and also,  by the Lebesgue Differentiation/Dominated Theorems, the
left hand side of this inequality tends to $0$ as $n \rightarrow
\infty$, hence:
\begin{equation}\label{insh2}
\mathcal{X}_{2}\rightarrow 
\int_{Q_T}F^{+}(u(x,t),v(x,t))g(x,t)\psi_{x}(x,t)dxdt\quad\mbox{as}\;\;n\rightarrow\infty.
\end{equation}
We now study the two terms $\mathcal{X}^{n}_{3}$ and
$\mathcal{X}^{n}_{4}$. From the fact that $g\in W^{1,\infty}(\bar{Q}_T)$, we remark
that for a.e. $(x,t,y,s)\in Q_{T}\times
Q_{T}$, we have:
\begin{equation}
g(x^{-},t^{-})-g(x^{+},t^{+})=
g_{x}(x^{-},t^{-})(-y/n)+g_{t}(x^{-},t^{-})(-s/n)+o\left(\frac{1}{n}\right).
\nonumber
\end{equation}
We also remark that the term $g_{x}(x^{+},t^{+})$ in
$\mathcal{X}^{n}_{4}$ could be replaced with $g_{x}(x^{-},t^{-})$, since
this adds a term that approaches $0$ as $n$ becomes large. This term will be omitted
throughout what follows and we denote the new $\mathcal{X}^{n}_{4}$ by
$\tilde{\mathcal{X}}^{n}_{4}$. From these two remarks, we rewrite
$\mathcal{X}^{n}_{3}$ and $\tilde{\mathcal{X}}^{n}_{4}$ to get:
\begin{equation}
\begin{aligned}
&\mathcal{X}^{n}_{3}=\int_{\mathcal{Q}_{4}}sgn^{+}(u(x^{+},t^{+})-v(x^{-},t^{-}))(f(u(x^{+},t^{+}))-f(v(x^{-},t^{-})))\\
&(yg_{x}(x^{-},t^{-})+sg_{t}(x^{-},t^{-}))\psi(x,t)\xi^{'}(y)\xi(s)dx\,dt\,dy\,ds+\mathcal{L}(n),
\end{aligned}
\end{equation}
where $\mathcal{L}(n)\rightarrow 0$ as $n\rightarrow \infty$,
and
\begin{equation}
\begin{aligned}
&\tilde{\mathcal{X}}^{n}_{4}=\int_{\mathcal{Q}_{4}}sgn^{+}(u(x^{+},t^{+})-v(x^{-},t^{-}))
(f(u(x^{+},t^{+}))-f(v(x^{-},t^{-})))\\
&\hspace{3cm}g_{x}(x^{-},t^{-})\psi(x,t)\xi(y)\xi(s)dx\,dt\,dy\,ds.
\end{aligned}
\end{equation}
The term $\mathcal{L}(n)$ will also be omitted for simplification and we
denote the new $\mathcal{X}^{n}_{3}$ by $\tilde{\mathcal{X}}^{n}_{3}$. Let
$\mathcal{X}^{n}_{34}=\tilde{\mathcal{X}}^{n}_{3}+\tilde{\mathcal{X}}^{n}_{4}$, hence:
\begin{equation}
\begin{aligned}
&\mathcal{X}^{n}_{34}=\overbrace{\int_{\mathcal{Q}_{4}}F^{+}(u(x^{+},t^{+}),v(x^{-},t^{-}))g_{x}(x^{-},t^{-})\psi(x,t)
(y\xi(y)\xi(s))_{y}dx\,dt\,dy\,ds}^{\mathcal{X}_{34}^{1n}}\\
&+\overbrace{\int_{\mathcal{Q}_{4}}F^{+}(u(x^{+},t^{+}),v(x^{-},t^{-})))g_{t}(x^{-},t^{-})
\psi(x,t)(s\xi(y)\xi(s))_{y}dx\,dt\,dy\,ds}^{\mathcal{X}_{34}^{2n}}.
\end{aligned}
\end{equation}
In $\mathcal{X}_{34}^{1n}$ and  $\mathcal{X}_{34}^{2n}$, the term $\psi(x,t)$ could be replaced with
$\psi(x^{-},t^{-})$, for this also adds a term getting small when $n\rightarrow
\infty$. We keep the same notations for $\mathcal{X}_{34}^{1n}$ and  $\mathcal{X}_{34}^{2n}$.
Since $y\xi(y)\xi(s)$ is a compactly supported smooth function
in $\mathcal{Q}_{4}$, we have:
\begin{equation}
\int_{\mathcal{Q}_{4}}F^{+}(u(x^{-},t^{-}),v(x^{-},t^{-}))g_{x}(x^{-},t^{-})\psi(x^{-},t^{-})
(y\xi(y)\xi(s))_{y}dx\,dt\,dy\,ds=0.
\end{equation}
Moreover, since $F^{+}(u,v)$ is Lipschitz continuous, we obtain:
\begin{equation}\label{giupbound}
\begin{aligned}
&\left|\mathcal{X}_{34}^{1n}-\int_{\mathcal{Q}_{4}}F^{+}(u(x^{-},t^{-}),v(x^{-},t^{-}))g_{x}(x^{-},t^{-})\psi(x^{-},t^{-})
(y\xi(y)\xi(s))_{y}dx\,dt\,dy\,ds\right|\\
&\leq
M_{f}(g_{x})^{\infty}\psi^{\infty}\int_{K_{\psi}}\int_{B_{1}^{2}}|u(x^{+},t^{+})-u(x^{-},t^{-})|dx\,dt\,dy\,ds,
\end{aligned}
\end{equation}
where $K_{\psi}$ is the support of $\psi$. Therefore, by the Lebesgue Differentiation/Dominated Theorems, we deduce that the
right hand side of (\ref{giupbound}) tends to $0$ as $n\rightarrow\infty$,
hence we have:
\begin{equation}\label{X341}
\mathcal{X}_{34}^{1n}\rightarrow 0\quad\mbox{as}\quad n\rightarrow\infty.
\end{equation}
In a similar way we can show that 
\begin{equation}\label{X342}
\mathcal{X}_{34}^{2n}\rightarrow 0\quad\mbox{as}\quad n\rightarrow\infty.
\end{equation}
From (\ref{insh1}), (\ref{insh2}), (\ref{X341}) and (\ref{X342}),
passing to the limit in (\ref{new_sum}) yields (\ref{relation_u_v}),
which concludes the proof of step 2.\\\\
{\bf Step 3: $u(x,t)\leq v(x,t)$ a.e. in $Q_T$.}\\\\
Let us first show that the function $A(t)$ defined in (\ref{A_t}) is
non-increasing a.e. in $(0,\min(T,\frac{a}{\omega}))$. Take $a>0$ and
recall that $\omega=g^{\infty}M_f$; let
$0<t_{1}<t_{2}<\min(T,\frac{a}{\omega})$,
$0<\eps<\min(t_{1},\min(T,\frac{a}{\omega}-t_{2})$, and
$\delta>0$. Consider the function $\phi\in C^{1}_{0}(\R_{+},[0,1])$ such
that $\phi(x)=1\;\forall x\in [0,a]$, $\phi(x)=0\;\forall x\in
[a+\delta,\infty)$, and $\phi^{'}<0$. Define $r_{\eps}$ by:
\begin{equation}\label{r_epsilon_t}
r_{\eps}(t)=\left\{
\begin{aligned}
&0 \quad&\mbox{if}&\quad  0\leq t\leq t_{1}-\eps\\
&\frac{t-(t_{1}-\eps)}{\eps} \quad &\mbox{if}&\quad t_{1}-\eps\leq t\leq 
t_{1}\\
&1   \quad&\mbox{if}&\quad t_{1}\leq t\leq t_{2}\\
&\frac{(t_{2}+\eps)-t}{\eps}  \quad&\mbox{if}& \quad t_{2}\leq t\leq 
t_{2}+\eps\\
&0  \quad&\mbox{if}&\quad t_{2}+\eps\leq t\leq \infty.
\end{aligned}
\right.
\end{equation}
One can take in (\ref{relation_u_v}) the permissible test function
$$\psi(x,t)=\phi(|x|+\omega t)r_{\eps}(t).$$
This yields:
\begin{equation}\label{E_EGH}
\begin{aligned}
&\overbrace{\frac{1}{\eps}\int_{t_{1}-\eps}^{t_{1}}\int_{\R}(u(x,t)-v(x,t))^{+}\phi(|x|+\omega
t)dxdt}^{E_{1}(\delta,\eps)}-\\
&\overbrace{\frac{1}{\eps}\int_{t_{2}}^{t_{2}+\eps}(u(x,t)-v(x,t))^{+}\phi(|x|+\omega
t)dxdt}^{E_{2}(\delta,\eps)}\geq E(\delta,\eps),
\end{aligned}
\end{equation}
with
\begin{equation}\label{E_eymgh}
\begin{aligned}
&E(\delta,\eps)&=&-\int_{0}^{T}\int_{\R}[\omega
(u(x,t)-v(x,t))^{+}+sgn^{+}((u(x,t)-v(x,t)))\times\\
&&(f&(u(x,t))-f(v(x,t)))\frac{x}{|x|}g(x,t)]\phi^{'}(|x|+\omega
t)r_{\eps}(t)dxdt.
\end{aligned}
\end{equation}
We claim that $E(\delta,\eps)\geq 0$. Indeed, since $\phi^{'}\leq 0$ and
$r_{\eps}\geq 0$, it suffices to
show that
\begin{equation}\label{1/2E_eymhg}
\begin{aligned}
&\omega(u(x,t)-v(x,t))^{+}+sgn^{+}((u(x,t)-v(x,t)))\times\\
&(f(u(x,t))-f(v(x,t)))\frac{x}{|x|}g(x,t)\geq
0\quad\mbox{a.e. in}\quad Q_T.
\end{aligned}
\end{equation}
Two cases can be considered, either $u(x,t)\leq v(x,t)$; in this case it
is easy to verify (\ref{1/2E_eymhg}), or $u(x,t)>v(x,t)$; in this
case we use, from the definition of $\omega$, the fact that
$$(f(u(x,t))-f(v(x,t)))\frac{x}{|x|}g(x,t)\geq -\omega
(u(x,t)-v(x,t)),$$
hence our claim holds. Relation (\ref{E_EGH}) now holds with 
$E(\delta,\eps)$
replaced with $0$. We regard the integrand term of $E_{1}(\delta,\eps)$ in 
(\ref{E_EGH}) and we notice that for $t_{1}-\eps<t<t_{1}$, we have:
$$(u(x,t)-v(x,t))^{+}\phi(|x|+\omega t)=(u(x,t)-v(x,t))^{+}\phi(|x|+\omega 
t)\mathbb{I}_{A_{\delta}},$$
where $\mathbb{I}_{A^{'}_{\delta}}$ is the characteristic function of the
set $A_{\delta}$ defined by:
$$A^{'}_{\delta}=\{(x,t);\; t_{1}-\eps<t<t_{1},\;0<|x|+\omega t<a+\delta\}.$$
Remark that the set $A^{'}_{\delta}$ shrinks, as $\delta$ becomes small, to
$$A^{'}=\{(x,t);\;t_{1}-\eps<t<t_{1},\;0<|x|+\omega t\leq a\}$$ with 
$\phi(|x|+\omega t)\equiv
1$ over $A$. It is easy now to see that as $\delta\rightarrow 0$
$$(u(x,t)-v(x,t))^{+}\phi(|x|+\omega
t)\mathbb{I}_{A^{'}_{\delta}}\rightarrow
(u(x,t)-v(x,t))^{+}\mathbb{I}_{A}\;\;\mbox{a.e. in}\;\;Q_T.$$
However, since $(u(x,t)-v(x,t))^{+}\in L^{\infty}(Q_T)$, we use the
Lebesgue Dominated Theorem to get:
\begin{equation}\label{E_1eps}
E_{1}(\delta,\eps)\rightarrow
\frac{1}{\eps}\int_{t_{1}-\eps}^{t_{1}}\int_{B_{a-\omega
    t}}(u(x,t)-v(x,t))^{+}dxdt\;\;\mbox{as}\;\;\delta\rightarrow 0,
\end{equation}
in other words,
\begin{equation}\label{E_1epsbis}
E_{1}(\delta,\eps)\rightarrow
\frac{1}{\eps}\int_{t_{1}-\eps}^{t_{1}}A(t)dt\;\;\mbox{as}\;\;\delta\rightarrow 
0,
\end{equation}
with $A(t)$ given by (\ref{A_t}). Similar arguments shows that:
\begin{equation}\label{E_2epsbis}
E_{2}(\delta,\eps)\rightarrow
\frac{1}{\eps}\int_{t_{2}}^{t_{2}-\eps}A(t)dt\;\;\mbox{as}\;\;\delta\rightarrow 
0.
\end{equation}
Note that $A\in L^{1}(0,T)$; let $t_{1}$ and $t_{2}$ be Lebesgue points
of the function $A$ such that $0<t_{1}<t_{2}<\min(T,\frac{a}{\omega})$,
one can easily deduce from (\ref{E_1epsbis}), (\ref{E_1epsbis}) and
(\ref{E_EGH}) letting $\eps$ tends to $0$ that
$$A(t_{1})\geq A(t_{2}),$$
hence $A$ is a.e. non-increasing. We use this property enjoyed by $A$ to
get the comparison principle. In fact, using the elementary identities:
\begin{equation}
\begin{aligned}
&(u-v)^{+}\leq (u-w)^{+}+(v-w)^{-}\\
&(u-v)^{-}\leq (u-w)^{-}+(v-w)^{+}
\end{aligned}
\nonumber
\end{equation}
$\forall\;u, v, w \in \R$, we calculate for a.e. $(x,t)\in
Q_T:$
$$(u(x,t)-v(x,t))^{+}\leq(u(x,t)-u^{0}(x))^{+}+(v(x,t)-v^{0}(x))^{-}+(u^{0}(x)-v^{0}(x))^{+}.$$
Since $u^{0}(x)\leq v^{0}(x)$ a.e. in $\R$, we get for a.e.
$(x,t)\in Q_T$:
\begin{equation}\label{l_a_s_t}
(u(x,t)-v(x,t))^{+}\leq(u(x,t)-u^{0}(x))^{+}+(v(x,t)-v^{0}(x))^{-}.
\end{equation}
Using (\ref{l_a_s_t}), for $\tau\in(0,T)$, we calculate:
\begin{equation}\label{conclusion}
\begin{aligned}
&\frac{1}{\tau}\int_{0}^{\tau}A(t)dt\leq\frac{1}{\tau}\int_{0}^{\tau}\int_{B_{a}}(u(x,t)-v(x,t))^{+}dxdt\leq\\
&\frac{1}{\tau}\int_{0}^{\tau}\int_{B_{a}}(u(x,t)-u^{0}(x))^{+}dxdt+
\frac{1}{\tau}\int_{0}^{\tau}\int_{B_{a}}(v(x,t)-v^{0}(x))^{-}dxdt.
\end{aligned}
\end{equation}
From (\ref{limfor_u0}), (\ref{limfor_v0}) and the passage to the
limit as $\tau\rightarrow 0$ in (\ref{conclusion}), we deduce that,
\begin{equation}\label{bet_l_eqn}
\frac{1}{\tau}\int_{0}^{\tau}A(t)dt\rightarrow 0 \quad\mbox{as}\quad
\tau \rightarrow 0.
\end{equation}
Thus, since $A$ is a.e. non-increasing on $(0,\tau)$, and $A(t)\geq
0$ for a.e. $t\in (0,\min(T,\frac{a}{\omega}))$, one then has
$$A(t)=0\quad\mbox{for a.e.}\quad t\in
\left(0,\min\left(T,\frac{a}{\omega}\right)\right).$$ 
Since $a$ is arbitrary, we deduce that,
$$u(x,t)\leq v(x,t)\quad\mbox{a.e. in}\quad Q_T.$$$\hfill{\Box}$
\begin{rem}\label{final_rem}
In \cite{egh2}, the entropy process solution $\mu(x,t,\alpha)$ was
proved to be independent of $\alpha$ for a divergence-free function
$g\in C^{1}(\bar{Q}_T)$. However, for the case of a general non
divergence-free function $g\in W^{1,\infty}(\bar{Q}_T)$, same result can
be shown by adapting the same proof as in \cite[Theorem 3]{egh2}
taking into account the slight modifications that could be deduced
from the proof of Theorem (\ref{entcomp}). More precisely, the
treatment of the two terms $\mathcal{X}^{n}_{3}$ and
$\mathcal{X}^{n}_{4}$ in Step 2.
\end{rem}
\end{section}
\noindent{\bf Acknowledgments}\\
The author would like to thank R. Monneau and C. Imbert for fruitful
discussions in the preparation of this paper. We also thank A.
El-Hajj, R. Eymard, N. Forcadel, M. Jazar and J. Vovelle for their
remarks. Finally, this work was partially supported by The
Mathematical Analysis and Applications Arab Network (MA3N) and by
The contract JC called "ACI
 jeunes chercheuses et jeunes chercheurs'' of the French Ministry of
 Research (2003-2007).

\vspace{1cm}
%\author{Hassan IBRAHIM~:}
\end{document}